\documentclass[twocolumn]{autartmod}
\usepackage[utf8]{inputenc}
\usepackage[T1]{fontenc}
\usepackage{mathpazo}
\usepackage{graphicx}
\usepackage{placeins}
\usepackage{amsmath,amssymb,bm}
\usepackage{xcolor}
\usepackage{silence}
\usepackage{subcaption}

\makeatletter
\renewcommand*\l@subsection{\@dottedtocline{2}{2.5em}{2.3em}} 
\makeatother

\allowdisplaybreaks 
\usepackage{ragged2e} 

\usepackage{tikz}
\usetikzlibrary{patterns} 
\usetikzlibrary{calc,positioning,shapes.geometric} 
\usepackage{pgfplots}
\pgfplotsset{compat=newest}
\usetikzlibrary{arrows.meta,pgfplots.fillbetween}
\newcommand{\coloredblock}[3][0.2cm]{\tikz[baseline=-0.75ex]{\node[fill=#2, draw=#3, minimum size=#1, inner sep=0pt] (n) {};}}

\usepackage{algorithm}
\usepackage{algpseudocode}

\colorlet{fadedred}{red!30!gray!80!white}
\colorlet{fadedblue}{blue!70!gray!80!white}
\colorlet{fadedgreen}{green!60!blue!80!gray}

\definecolor{beige}{RGB}{245,245,220}

\usetikzlibrary{shadows}
\usepackage[many]{tcolorbox}
\newtcolorbox{eqbox}{
  colback=beige!30,     
  colframe=black!30,    
  boxrule=0.7pt,        
  arc=5pt,              
  boxsep=0pt,           
  left=0pt, right=5pt,  
  top=-6pt, bottom=-10pt
}
\newtcolorbox{eqboxb}{
  colback=beige!30,     
  colframe=black!30,    
  boxrule=0.7pt,        
  arc=5pt,              
  boxsep=0pt,           
  left=0pt, right=5pt,  
  top=-6pt, bottom=5pt
}

\usepackage{hyperref}
\hypersetup{
    colorlinks=true,
    linkcolor=fadedblue, 
    citecolor=fadedgreen,
    urlcolor=fadedgreen,
    pdftitle={DT-STL}
}

\renewenvironment{pf}{\par\textcolor{gray}{\textit{\textbf{Proof:}~}}}{\hfill\coloredblock{black!50}{black!50}\par}

\renewcommand{\dot}[1]{\overset{\text{{\large{.}}}}{#1}}

\usepackage{comment}
\usepackage{todonotes}
\usepackage{cite}
\usepackage{textcomp}

\usepackage{empheq}
\usepackage{float}

\newcommand{\bYY}[1]{\color{black} #1 \color{black}}

\newcommand{\fbs}[1]{\color{black} #1 \color{black}}

\usepackage{graphicx}

\newtheorem{lemma}{Lemma}
\newtheorem{definition}{Definition}
\newtheorem{assumption}{Assumption}

\newtheorem{proof}{Proof}
\newtheorem{example}{Example}

\usepackage{array}

\newcommand{\greentick}{\textcolor{green!70!black}{\checkmark}}
\newcommand{\yellowcircle}{\textcolor{yellow!80!black}{$\circ$}}
\newcommand{\redcross}{\textcolor{red}{$\times$}}

\newcommand{\intv}[2]{[#1{:}#2]}

\newcommand{\Hsquare}{%
  \text{\fboxsep=-.2pt\fbox{\rule{0pt}{1ex}\rule{1ex}{0pt}}}%
}
\edef\endfrontmatter{%
  \unexpanded\expandafter{\endfrontmatter}
  \noexpand\endNoHyper 
}
\begin{document}
\begin{frontmatter}
\runtitle{D-GSAR}
\title{%
Optimization with Temporal and Logical Specifications \\
via Generalized Mean-based Smooth Robustness Measures
\!\!\!\!
\thanksref{footnoteinfo}
}
\thanks[footnoteinfo]{This research was supported by AFOSR grant FA9550-20-1-0053 and ONR grant N00014-20-1-2288; Government sponsorship is acknowledged. \\ Corresponding author: Samet Uzun.}
\author{Samet Uzun}$^{*}$\ead{samet@uw.edu},~~%
\author{Purnanand Elango}$^{*}$\ead{pelango@uw.edu},~~%
\author{Pierre-Lo{\"i}c Garoche}$^{\dagger}$\ead{pierre-loic.garoche@enac.fr},~~%
\author{Beh{\c{c}}et A{\c{c}}{\i}kme{\c{s}}e}$^{*}$\ead{behcet@uw.edu}%
\address{$^{*}$William E.\ Boeing Department of Aeronautics \& Astronautics, University of Washington, Seattle, WA, 98195}
\address{$^{\dagger}${\'E}cole Nationale de l'Aviation Civile, Universit{\'e} de Toulouse, France}
\begin{abstract}
    This paper introduces a generalized mean-based $\mathcal{C}^1$-smooth robustness measure over discrete-time signals (D-GMSR) for signal temporal logic (STL) specifications.
    In conjunction with its $\mathcal{C}^1$-smoothness, D-GMSR is proven to be both {\em sound and complete}.
    Furthermore, it demonstrates favorable gradient properties and addresses {\em locality and masking} problems, which are critical for numerical optimization.
    The $\mathcal{C}^1$-smoothness of the proposed formulations enables the implementation of robust and efficient numerical optimization algorithms to solve problems with STL specifications while preserving their theoretical guarantees.
    The practical utility of the proposed robustness measure is demonstrated on two real-world trajectory optimization problems: i) quadrotor flight, and ii) autonomous rocket landing.
    A sequential convex programming (SCP) framework, incorporating a convergence-guaranteed optimization algorithm (the prox-linear method) is used to solve inherently non-convex trajectory optimization problems with STL specifications.
    The implementation is available at \url{https://github.com/UW-ACL/D-GMSR}
    \vspace{-0.3cm}
\end{abstract}

\end{frontmatter}
\section{Introduction} \label{sec:intro}
\vspace{-0.15cm}
Temporal logic specifications \cite{baier2008principles} provide a systematic approach for encoding high-level specifications and complex system behaviors in mission planning and control.
The application of temporal logic specifications in dynamical systems can be found in autonomous vehicles \cite{sahin2020autonomous}, high-dimensional manipulators \cite{kurtz2020trajectory}, and multi-agent systems \cite{djeumou2022probabilistic}, including mobile robots \cite{sahin2019multirobot,buyukkocak2021planning}, and UAVs \cite{pant2018fly,bacspinar2019mission}.
Temporal logic specifications are also integrated into learning algorithms to ensure that results satisfy certain properties in prediction \cite{ma2020stlnet} and policy synthesis \cite{aksaray2016q,alshiekh2018safe}. 
Notably, useful formulations have emerged to facilitate systematic incorporation of such constraints \cite{leung2023backpropagation}.

\vspace{-0.15cm}
Linear temporal logic (LTL) is one of the most popular temporal logic frameworks for expressing design requirements for discrete event systems \cite{pnueli1977temporal}.
Finite state automata \cite{alur1994theory} have been traditionally employed to control discrete-time dynamical systems with LTL specifications.
A receding horizon framework is proposed in \cite{wongpiromsarn2012receding} to alleviate the computational complexity of synthesizing finite state automata, along with the associated software toolbox \cite{wongpiromsarn2011tulip}.
\cite{karaman2008optimal} introduced a general technique, formulating LTL specifications as mixed-integer linear constraints, enabling optimization algorithms to address these problems.
Metric temporal logic \cite{koymans1990specifying} and metric interval temporal logic \cite{alur1996benefits} serve as natural extensions of LTL, to extend the discrete-time specifications to continuous time.   
Signal temporal logic (STL) combines the dense temporal modalities of metric interval temporal logic (MITL) with predicates specifying properties over real-valued signals, motivated by the application of verifying and monitoring temporal properties of continuous-time signals in continuous and hybrid systems \cite{maler2004monitoring}. A robust semantic for MITL, which generates a real number to quantify the degree of specification satisfaction instead of a Boolean evaluation, is introduced in \cite{fainekos2009robustness}.
A similar robustness measure, called space robustness (SR), and a time robustness concept aimed at improving robustness against timing uncertainties within STL is introduced in \cite{donze2010robust}.
SR utilizes $\min$ and $\max$ functions, \fbs{which we refer to as key functions of SR,} for logical and temporal operators of STL.
Hence, the introduced robustness measures allow the encoding of temporal specifications as cost functions, transforming the synthesis problem into an optimization problem \cite{lindemann2018control,belta2019formal}.
\bYY{While SR is initially formulated for continuous-time signals \cite{donze2010robust}, a discrete-time counterpart (D-SR), which resembles the robustness measure proposed for MITL in \cite{fainekos2009robustness}, has been explored in related literature \cite{raman2014model,sadraddini2015robust,rodionova2021time,rodionova2022combined,pant2017smooth,haghighi2019control,gilpin2020smooth,mao2022successive,lindemann2019robust,mehdipour2019arithmetic,mehdipour2020specifying}.}
The inherent design of the introduced robustness measures leads to solutions in the form of mixed-integer programming (MIP) problems for both D-SR \cite{raman2014model,sadraddini2015robust} and time robustness \cite{rodionova2021time,rodionova2022combined}.
The formulation of D-SR with the $\min$ and $\max$ functions has been addressed by smoothing these functions using the well-known log-sum-exponential (LSE) functions in \cite{pant2017smooth}. The work is further extended in \cite{haghighi2019control} to address the notion of cumulative robustness.
The LSE approximation is asymptotically sound and asymptotically complete, serving as an under-approximation for the $\min$ function and an over-approximation for the $\max$ function.
The LSE approximation for the $\max$ function is replaced with an under-approximation in \cite{gilpin2020smooth}, resulting in a sound, asymptotically complete and smooth space robustness (D-SSR) for STL.
A polynomial smooth approximation is used in \cite{mao2022successive} for a more accurate approximation of the $\min$ and $\max$ functions. 
\fbs{Although the number of variables of these functions that approximates the $\min$ and $\max$ functions cannot be greater than two,} recursive modeling of STL constraints in this paper causes all $\min$ and $\max$ functions to perform on only two variables.
While the robust semantics introduced for STL specifications demonstrate utility in verification and monitoring \cite{donze2010robust}, the smooth versions may pose challenges for optimization, primarily related to concerns about locality and masking issues defined in \cite{mehdipour2019average}.
Locality is defined as the dependency of the robustness measure on signals at a single time instant, whereas masking describes the situation where the robustness measure is solely determined by one of the subformulas combined with  \textit{conjunction} and  \textit{disjunction} operators.
Utilizing an averaging approach that considers the entire signal and subformulas, instead of relying on specific time instances or subformulas through the $\min$ and $\max$ functions, proves beneficial in mitigating these issues.
Averaged-STL is introduced in \cite{akazaki2015time} to address the falsification problem in continuous-time signals. This approach involves creating a weighted time average over specific STL operators. However, it does not eliminate the use of the $\min$ and $\max$ functions, 
resulting in the continuation of non-smooth STL semantics.
In \cite{lindemann2019robust}, discrete average space robustness (DASR) and its simplified version (DSASR) are introduced as computationally more tractable alternatives to D-SR. However, both DASR and DSASR incorporate $\min$ and $\max$ functions, leading to non-smoothness in certain STL operators. Additionally, ensuring the soundness of the STL formula necessitates imposing additional constraints.
The arithmetic-geometric mean robustness (AGM), presented in \cite{mehdipour2019arithmetic}, emphasizes both the frequency of satisfaction of STL specifications and the robustness of each specification at every time instant.
The AGM also addresses the locality and masking problems. However, the introduced robustness measure switches between arithmetic and geometric means with if-else conditions to satisfy the STL specification, leading to non-smoothness. Additionally, the robustness degree can only cover a specific set of STL specifications.
Both D-SR and AGM are generalized in \cite{mehdipour2020specifying} by introducing weights associated with STL operators to enable the definition of user preferences for priorities and importance. 

\vspace{-0.15cm}
This paper introduces a generalized mean-based $\mathcal{C}^1$-smooth robustness measure over discrete-time signals (D-GMSR). D-GMSR is constructed by replacing the key functions of D-SR, $\min$, and $\max$ functions, with two $\mathcal{C}^1$-smooth functions that are specific compositions of generalized means. These functions have the same sign with $\min$ and $\max$ functions, respectively. 
Therefore, D-GMSR is both sound and complete, making it reliable and non-conservative.
Furthermore, the first-order information of these functions satisfies some desirable properties critical for numerical optimization purposes and addresses the issues of locality and masking.
Table \ref{tab:evaluation} compares D-GMSR with existing robustness measures from the literature.
The $\mathcal{C}^1$-smoothness of the robustness measure facilitates the application of a broad range of numerical solvers for optimization problems with STL specifications while preserving their theoretical guarantees.
We provide some illustrative numerical examples of D-GMSR and demonstrate its practical utility in trajectory optimization for quadrotor flight and autonomous rocket landing problems solved via sequential convex programming (SCP) framework \cite{mao2016successive,malyuta2021advances,malyuta2022convex,ctcs2024}. 
The implementation is available at \url{https://github.com/UW-ACL/D-GMSR}
\begin{table}[htbp]
  \centering
  \renewcommand{\arraystretch}{1.2}
  \resizebox{\columnwidth}{!}{%
  \begin{tabular}
  {|c|c|c|c|c|c|}
    \hline
    & \cite{donze2010robust}  & \cite{pant2017smooth} & \cite{gilpin2020smooth}  & \cite{mehdipour2019arithmetic}  & D-GMSR \\
    \hline
    Smoothness          & \redcross  & \greentick    & \greentick    & \redcross  & \greentick \\
    \hline
    Soundness           & \greentick & \yellowcircle     & \greentick    & \greentick & \greentick \\
    \hline
    Completeness        & \greentick & \yellowcircle & \yellowcircle & \greentick & \greentick \\
    \hline
    Locality \& Masking & \redcross  & \redcross     & \redcross     & \greentick & \greentick \\
    \hline
  \end{tabular}
  }
  \caption{Comparison with the previous robustness measures \\ (\yellowcircle: Satisfied Asymptotically)}
  \label{tab:evaluation}
\end{table}

\vspace{-0.35cm}
The paper is organized as follows. 
Section \ref{sec:stl} introduces the syntax of STL and smooth space robustness, quantifying STL specification satisfaction, yet highlighting its limitations resulting in undesirable optimization outcomes where the specification is not met.
These limitations are addressed by the proposed generalized mean-based $\mathcal{C}^1$-smooth robustness measure (D-GMSR) in Section \ref{sec:gsar}. 
Section \ref{sec:num_sim} demonstrates the practical utility of D-GMSR in trajectory optimization for quadrotor flight and rocket landing problems solved within an SCP framework.
The paper is concluded in Section \ref{sec:conc}.

\textit{Notation:} 
The real numbers are denoted as $\mathbb{R}$, where $\mathbb{R}_+$ and $\mathbb{R}_{++}$ represent the sets of non-negative and positive real numbers, respectively. 
Similarly, $\mathbb{Z}_+$ and $\mathbb{Z}_{++}$ represent the sets of non-negative and positive integers. 
$\mathbb{R}^{n \times m}$ and $\mathbb{R}^n$ denote the sets of real $n \times m$ matrices and $n \times 1$ vectors, respectively. 
$\mathrm{diag}(A) \in \mathbb{R}^n$ represents the diagonal elements of a matrix $A \in \mathbb{R}^{n \times n}$.
$\bm{0}$ denotes a zero vector, while $\bm{1}$ represents a vector of ones in appropriate dimensions.  
The vector $e_i$ is defined such that its $i$th element is $1$ and all other elements are $0$.
The cardinality of set $S$ is denoted by $\mathrm{card}(S)$. 
$\top$ and $\bot$ represent boolean true and false, respectively.
The space of continuous and $k$-times differentiable functions on a domain $X$ is denoted by $\mathcal{C}^k(X)$, with members referred to as $\mathcal{C}^k$-smooth.
A vector $x \in \mathbb{R}^n$ is represented as $(x_1, x_2, \dots , x_n)$.
$\intv{a}{b} := \{ k \in \mathbb{Z}_{+} : a \leq k \leq b; a,b \in \mathbb{Z}_{+} \}$. 
$x_{\intv{a}{b}} := (x_a, x_{a+1}, \dots, x_b)$.
The functions $\min(x)$ and $\max(x)$ return the minimum and maximum elements of the $x \in \mathbb{R}^n$ vector, respectively.
$|x|_{\Hsquare}$ is defined as either $|x|_{+}$ or $|x|_{-}$. 
$|x|_{\Hsquare}^p = (f((x_1, 0))^p, \dots, f((x_n, 0))^p)$, where $f = \max$ if $\Hsquare = +$, and $\min$ otherwise.
\vspace{-0.2cm}
\section{Preliminaries} \label{sec:stl}
\vspace{-0.2cm}
In this section, we provide a concise overview of the syntax and a robust semantic of STL tailored for discrete-time signals.
\vspace{-0.15cm}
\subsection{Signal Temporal Logic (STL)}
\vspace{-0.15cm}
Let $x: \mathbb{Z}_+ \to \mathbb{R}^n$ be a discrete-time signal and $\mu$ be a 
Boolean-valued predicate such that
\begin{align*}
    \mu :=& \; (f(x) \geq 0) \\ 
    =&
    \begin{cases}
        \top & \text{if $f(x) \geq 0$} \\
        \bot & \text{otherwise}.
    \end{cases} \nonumber
\end{align*}
where $f: \mathbb{R}^n \to \mathbb{R}$ is a real-valued function referred to as the predicate function. The following assumption 
is needed for numerical optimization purposes.
\vspace{-0.15cm}
\begin{assumption} \label{asm:pred_smth}
    All predicate functions are $\mathcal{C}^1$-smooth.
\end{assumption}
\vspace{-0.15cm}
Let $I$ be a non-empty time interval such that $I = \intv{a}{b}$. 
The syntax of STL is defined recursively as follows:
\begin{equation} \label{eq:stl_syntax}
    \varphi ::= \mu \; | \; \neg \varphi_1 \; | \; \varphi_1 \wedge \varphi_2 \; | \; \varphi_1 \bm{U}_I \varphi_2
\end{equation}
where \bYY{$\varphi$ is the STL specification (or STL formula), $\varphi_1$ and $\varphi_2$ are the STL sub-specifications (or STL sub-formulas)}, $\neg$ and $\wedge$ are the Boolean \textit{negation} and \textit{conjunction}, respectively, and $\bm{U}_I$ is the \textit{until} temporal operator over bounded interval $I$. \bYY{In order to construct an STL specification, different predicates or STL sub-formulas are recursively combined using the Boolean \textit{negation}, \textit{conjunction}, and \textit{until} operators as in the examples in Table \ref{tab:stl_opt}.}
\begin{table}[H]
    \centering
    \renewcommand{\arraystretch}{1.1}
    \begin{tabular}{|cc|c|}
    \hline
    \multicolumn{2}{|c|}{STL specification}     & Formula \\ \hline
    \multicolumn{1}{|c|}{\textit{disjunction}} & $\varphi_1 \vee \varphi_2$ & $\neg ( \neg \varphi_1 \wedge \neg \varphi_2 )$ \\ \hline
    \multicolumn{1}{|c|}{\textit{implication}} & $\varphi_1 \implies \varphi_1$ & $\neg \varphi_1 \vee \varphi_2$ \\ \hline
    \multicolumn{1}{|c|}{\textit{eventually}} & $\bm{F}_I \varphi$ & $\top \bm{U}_I \varphi$ \\ \hline
    \multicolumn{1}{|c|}{\textit{always}} & $\bm{G}_I \varphi$ & $\neg \bm{F}_I \neg \varphi$ \\ \hline
    \end{tabular}
    \caption{Construction of some STL specifications}
    \label{tab:stl_opt}
\end{table}
\vspace{-0.5cm}
Formally, the semantics of an STL formula $\varphi$ defines what it means for a signal $x$ to satisfy $\varphi$ at time point $k$, denoted as $(x,k) \vDash \varphi$.
Conversely, if the specification is not satisfied, it is denoted as $(x,k) \nvDash \varphi$.
\begin{rem}
The $\bm{U}_I$ operator in the STL syntax (\ref{eq:stl_syntax}), and consequently the $\bm{F}_I$ and $\bm{G}_I$ operators, can be represented over discrete-time signals by utilizing the first three elements, namely the predicate, \textit{negation}, and \textit{conjunction}.
\end{rem}

\vspace{-0.15cm}
\subsection{Space Robustness for STL}
\vspace{-0.15cm}
While determining if the signal meets the STL specification with a Boolean evaluation, the robust semantics yield a real number that quantifies the degree of satisfaction of a specification. 
A space robustness (SR) concept, a robustness measure for STL, is presented in \cite{donze2010robust} by utilizing $\min$ and $\max$ functions for logical and temporal operators of STL. In this approach, temporal specifications can be represented as cost functions, transforming a synthesis problem into an optimization problem \cite{lindemann2018control, belta2019formal}.
While SR is initially developed for continuous-time signals, its discrete-time version (D-SR) has been explored in various studies \cite{raman2014model, sadraddini2015robust, rodionova2021time, rodionova2022combined, pant2017smooth, haghighi2019control, gilpin2020smooth, mao2022successive, lindemann2019robust, mehdipour2019arithmetic, mehdipour2020specifying}.
\begin{definition} \label{def:dtsr} 
\fbs{(D-SR \cite{donze2010robust, raman2014model})} The space robustness \cite{donze2010robust} defines a real-valued function $\rho^{\varphi}$ of signal $x$ and time $k$ such that $(x,k) \models \varphi \iff \rho^\varphi (x, k) \geq 0$. 
The complete D-SR is defined as follows:
\begin{align*}
    &\rho^{\mu} (x,k) := f(x_k)\\
    &\rho^{\neg \varphi} (x,k) := - \rho^\varphi (x,k)\\
    &\rho^{\varphi_1 \wedge \varphi_2} (x,k) := \min (( \rho^{\varphi_1} (x,k), \rho^{\varphi_2} (x,k) ))\\
    &\rho^{\varphi_1 \vee \varphi_2} (x,k) := \max (( \rho^{\varphi_1} (x,k), \rho^{\varphi_2} (x,k) ))\\
    &\rho^{\varphi_1 \implies \varphi_2} (x,k) := \max (( -\rho^{\varphi_1} (x,k), \rho^{\varphi_2} (x,k) ))\\
    & \rho^{ \bm{F}_{[a : b]} \varphi}  (x,k) := \max(( \rho^{ \varphi }(x, k+a), \rho^{ \varphi }(x, k+a+1), \\
    &\hspace{3.3cm}\dots, \rho^{ \varphi }(x, k+b) ))\\
    & \rho^{ \bm{G}_{[a : b]} \varphi}  (x,k) := \min(( 
        \rho^{ \varphi }(x, k+a), \rho^{ \varphi }(x, k+a+1),\\
    &\hspace{3.3cm}\dots, \rho^{ \varphi }(x, k+b)))\\
    & \rho^{ \varphi_1 \bm{U}_{[a : b]} \varphi_2}  (x,k) := \max(( \min(( \rho^{ \bm{G}_{[a : a]} \; \varphi_1}(x, k), \\
    &\hspace{4.76cm}\rho^{ \varphi_2 }(x, k+a) )), \\
    &\hspace{3.84cm}\min(( \rho^{ \bm{G}_{[a : a+1]} \; \varphi_1}(x, k), \\
    &\hspace{4.76cm}\rho^{ \varphi_2 }(x, k+a+1) )), \\
    &\hspace{4.175cm}\vdots \\
    &\hspace{3.84cm}\min(( \rho^{ \bm{G}_{[a : b]} \; \varphi_1}(x, k), \\
    &\hspace{4.76cm}\rho^{ \varphi_2 }(x, k+b) )) ))
\end{align*}
\end{definition}

The $\min$ and $\max$ functions are key functions of D-SR, and any STL specification can be formulated by composing the predicate functions with these key functions.
However, the non-smoothness of the $\min$ and $\max$ functions results in optimization problems involving STL specifications becoming MIP problems \cite{raman2014model, sadraddini2015robust}.
In \cite{gilpin2020smooth}, a smooth version of D-SR (D-SSR) is proposed by replacing $\min$ and $\max$ functions with their under-approximations as follows:
\begin{subequations} \label{eq:app_minmax}
\begin{align}
    \fbs{\widetilde{\min}_{\kappa}}(x) &:= - \frac{1}{\kappa} \mathrm{log} \bigg( \sum_{i=1}^n e^{-\kappa x_i} \bigg) \\
    \fbs{\widetilde{\max}_{\kappa}}(x) &:= \frac{\sum_{i=1}^n x_i e^{\kappa x_i}}{ \sum_{i=1}^n e^{\kappa x_i} }
\end{align}
\end{subequations}
where $x \in \mathbb{R}^n$ and $\kappa >0$ is a sharpness parameter determining the accuracy of the approximation. 
Thus, a smooth, sound, and asymptotically complete robustness measure $\tilde{\rho}^{\varphi}(x,k)$ is achieved for STL, enabling the utilization of gradient-based algorithms to solve optimization problems involving STL specifications.

\subsection{Gradient properties of the D-SSR} \label{sec:grad_dssr}
In D-SR, the values of multiple predicate functions at various time instances can be composed using the $\min$ and $\max$ functions, such as
\begin{equation*}
    \max((\min ((f_1(x_1), f_2(x_1))), f_1(x_2))).
\end{equation*}
Due to the nature of $\min$ and $\max$ functions, the robustness value of an STL specification is typically determined by the value of a specific predicate function at a single time instant. The values of all other predicate functions and values of predicate functions at other time instances are ignored, which are called {\em locality and masking}, respectively, as introduced in \cite{mehdipour2019average}.
Consequently, the derivatives of  $\fbs{\widetilde{\min}_{\kappa}}$ and $\fbs{\widetilde{\max}_{\kappa}}$ functions with respect to one of the input variables are nearly zero if  the minimum or maximum value are not attained by this variable among others. 

To demonstrate this, suppose that $f_i: \mathbb{R}^n \to \mathbb{R}$ is a function of $x \in \mathbb{R}^n$ for all $i \in \intv{1}{m}$, $S_{\textrm{min}} = \{ i : f_i(x) \leq f_j(x), \forall j \in \intv{1}{m} \} $, and $S_{\textrm{max}} = \{ i : f_i(x) \geq f_j(x), \forall j \in \intv{1}{m} \} $. For sufficiently large $\kappa$, if $y = (f_1(x), f_2(x), \dots, f_m(x) ) $, then 
\begin{align*} \label{eq:grad_minmax}
    \nabla_x \fbs{\widetilde{\min}_{\kappa}}( y ) &= \sum_{i=1}^m a_i \nabla_x f_i(x) \text{ and}\\
    \nabla_x \fbs{\widetilde{\max}_{\kappa}}( y ) &= \sum_{i=1}^m b_i \nabla_x f_i(x), \text{ where}   
\end{align*}
\begin{align*}
    a_i &\approx
    \begin{cases}
        \frac{1}{ \mathrm{card}(S_{\textrm{min}}) } & \text{if } i \in S_{\textrm{min}} \\
        0 & \text{otherwise}
    \end{cases}, \\
    b_i &\approx
    \begin{cases}
        \frac{1}{ \mathrm{card}(S_{\textrm{max}}) } & \text{if } i \in S_{\textrm{max}} \\
        0 & \text{otherwise}
    \end{cases}.
\end{align*}
Therefore, using D-SSR  to account for an STL specification within a gradient-based iterative optimization algorithm leads to the consideration of only the gradient of one predicate function at a single time instant, $f_i(x_k)$.

\begin{figure}
    \centerline{\includegraphics[scale=0.43]{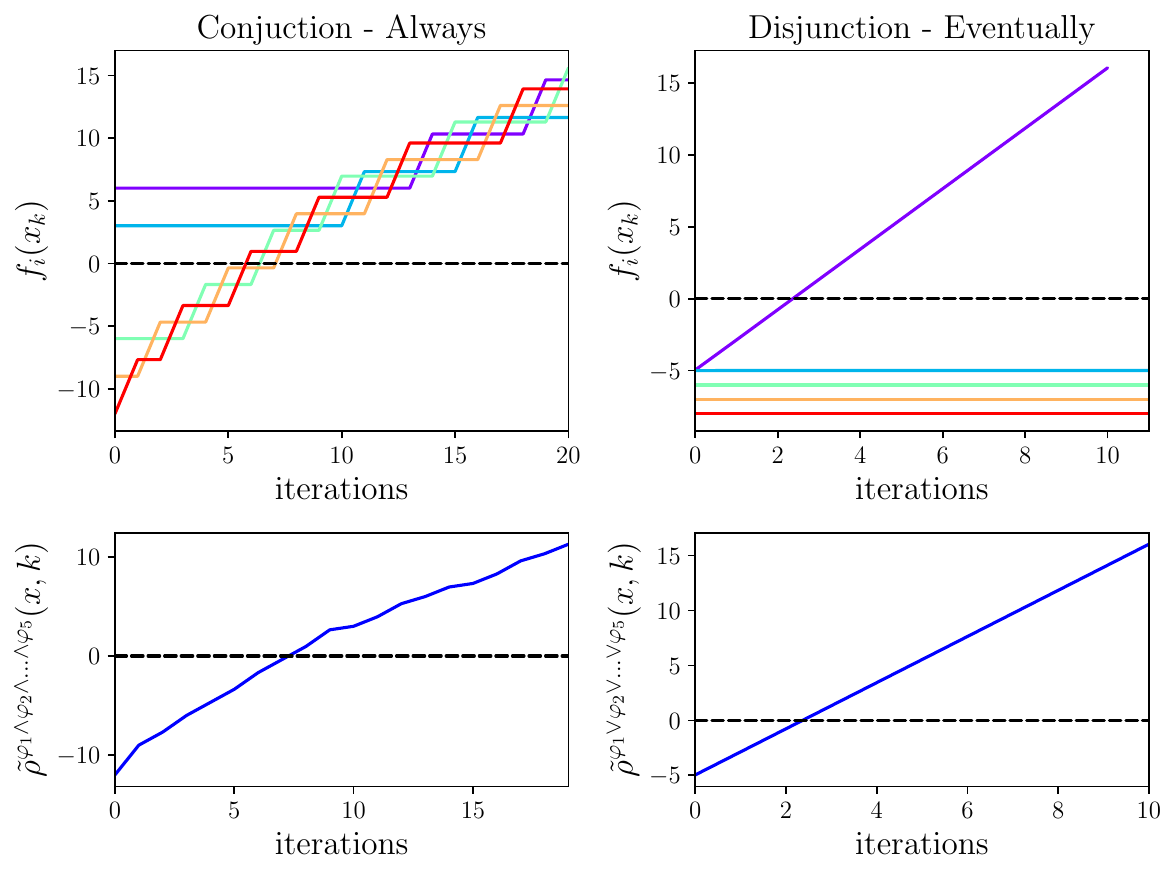}}
    \caption{
    The changes in the values of $f_i(x_k)$ (or $f(x_i)$) and the corresponding robustness values while maximizing the \textit{conjunction} (or \textit{always}) and \textit{disjunction} (or \textit{eventually}) operators using D-SSR 
    are presented on the left and right sides of the figure, respectively.
    \bYY{The color gradient, ranging from purple to red in the figures above, represents the values of the functions from $f_1(x_k)$ to $f_5(x_k)$ (or $f(x_1)$ to $f(x_5)$), respectively.}
    }
    \label{fig:dssr}
\end{figure}

\begin{example}
    In Figure \ref{fig:dssr}, we illustrate the variations in the predicate function values during the maximization of $\tilde{\rho}^{\varphi_1 \wedge \varphi_2 \wedge \dots \wedge \varphi_5} (x,k)$ and $\tilde{\rho}^{\varphi_1 \vee \varphi_2 \vee \dots \vee \varphi_5} (x,k)$ functions, respectively. \fbs{The sharpness parameter determining the accuracy of the approximation of the $\min$ and $\max$ functions \eqref{eq:app_minmax} is set to $\kappa=25$.} Here, $\tilde{\rho}^{\varphi_i}(x,k) = f_i(x_k) \in \mathbb{R}$ and $f_i : \mathbb{R} \to \mathbb{R}$ for all $i \in \intv{1}{5}$. Therefore,
    \begin{align*}
        \tilde{\rho}^{\varphi_1 \wedge \varphi_2 \wedge \dots \wedge \varphi_5} (x,k) &= \fbs{\widetilde{\min}_{\kappa}} ( y ), \\
        \tilde{\rho}^{\varphi_1 \vee \varphi_2 \vee \dots \vee \varphi_5} (x,k) &= \fbs{\widetilde{\max}_{\kappa}} ( y ),
    \end{align*}
    where $y=( f_1(x_k), f_2(x_k), \dots, f_5(x_k)  )$.
    Alternatively, maximizing $\tilde{\rho}^{ \bm{G}_{[1 : 5]} \varphi} (x,0)$ and $\tilde{\rho}^{ \bm{F}_{[1 : 5]} \varphi} (x,0)$ functions, respectively, where $\tilde{\rho}^{\varphi}(x,k) = f(x_k) \in \mathbb{R}$ for all $k \in \intv{1}{5}$ and $f : \mathbb{R} \to \mathbb{R}$, leads to similar changes in predicate function values across different time instances since
    \begin{align*}
        \tilde{\rho}^{\bm{G}_{[1 : 5]} \varphi} (x,0) &= \fbs{\widetilde{\min}_{\kappa}} ( y ), \\
        \tilde{\rho}^{\bm{F}_{[1 : 5]} \varphi} (x,0) &= \fbs{\widetilde{\max}_{\kappa}} ( y ),
    \end{align*}
    where $y=( f(x_1), f(x_2), \dots, f(x_5)  )$.
    
    During the maximization of $\tilde{\rho}^{\varphi_1 \wedge \varphi_2 \wedge \dots \wedge \varphi_5} (x,k)$ function, only the value of the predicate function that attains the minimum value increases notably with each iteration and the function with the minimum value changes after a certain number of iterations.
    Hence, the values of all predicate functions become positive after a large number of iterations. On the other hand, during the maximization of $\tilde{\rho}^{\varphi_1 \vee \varphi_2 \vee \dots \vee \varphi_5} (x,k)$ function, only the value of the predicate function that attains the maximum value increases notably with each iteration. Therefore, the predicate function with the highest value does not change and the relevant function becomes positive after a large number of iterations.
    These characteristics of $\fbs{\widetilde{\min}_{\kappa}}$ and $\fbs{\widetilde{\max}_{\kappa}}$ functions may lead to undesirable results, i.e.,  the STL specification is not met, as illustrated by an example in Section \ref{sec:num_sim}. 
\end{example}

These undesirable characteristics might be alleviated by introducing a small smoothing parameter $\kappa$ for the $\min$ and $\max$ functions in the initial iterations of the numerical optimization algorithm, gradually increasing $\kappa$ with subsequent iterations to achieve a more accurate approximation similar to the homotopy approach employed in the SCP algorithm \cite{malyuta2023fast}.
However, this approach requires the inclusion of additional decision variables in the algorithm and may potentially disrupt its useful convergence guarantees.
\section{
Generalized Mean-based Smooth Robustness over Discrete-time Signals (D-GMSR)
} \label{sec:gsar}
To remedy the issues in previously proposed robustness measures, this section introduces generalized mean-based $\mathcal{C}^1$-smooth robustness measure over discrete-time signals (D-GMSR).
We first define functions that are the key for D-GMSR and then introduce D-GMSR by replacing $\min$ and $\max$ functions of D-SR in Definition \ref{def:dtsr} with the proposed key functions.
Here we also provide proofs of the $\mathcal{C}^1$-smoothness, soundness, and completeness properties of the D-GMSR.
We then present an analysis of the gradients of the proposed key functions and provide 
some numerical illustrative examples for the application of D-GMSR.

We propose the following $\mathcal{C}^1$-smooth function, constructed through generalized means, characterized by the same sign as the $\min$ function. This proposed function will replace the $\min$ function of D-SR in Definition \ref{def:dtsr}.
\begin{equation} \label{eq:dgsar_key}
    {}^{\wedge} h_{p, w}^{\epsilon}(y) := 
    \bigg( M_{0, w}^{\epsilon}(|y|_{+}^2) \bigg)^{\frac{1}{2}}
    \hspace{-0.2cm} - \bigg( M_{p, w}^{\epsilon}(|y|_{-}^2) \bigg)^{\frac{1}{2}},
\end{equation}
where
\begin{subequations} 
\begin{align} 
    M_{0, w}^{\epsilon}(z) &:= \Big( \epsilon^{\bm{1}^T w} + \prod_{i=1}^n z_i^{w_i} \Big)^{1/\bm{1}^T w} \\
    M_{p, w}^{\epsilon}(z) &:= \Big( \epsilon^p + \frac{1}{\bm{1}^T w} \sum_{i=1}^n w_i z_i^p \Big)^{1/p}
\end{align}
\label{eq:key_fcn}%
\end{subequations}
$y \in \mathbb{R}^n$, $z \in \mathbb{R}^n_+$, $\epsilon \in \mathbb{R}_{++}$, $p \in \mathbb{Z}_{++}$, $w \in \mathbb{Z}_{++}^n$. We define another key function that will replace the $\max$ function of D-SR in Definition \ref{def:dtsr} as ${}^{\vee} h_{p, w}^{\epsilon}(y) := -{}^{\wedge} h_{p, w}^{\epsilon}(-y)$.
\begin{rem} \label{rem:key}
    If $\epsilon = 0$, then $M_{0, w}^{\epsilon}(z)$ and $M_{p, w}^{\epsilon}(z)$ are weighted geometric mean and weighted power mean of $z$, respectively. Furthermore, $M_{0, w}^{\epsilon}, M_{p, w}^{\epsilon} : \mathbb{R}^n_+ \to [\epsilon, \infty)$, and $M_{0, w}^{\epsilon}, M_{p, w}^{\epsilon} \in \mathcal{C}^{\infty} (\mathbb{R}_+^n) $. The values of the ${}^{\wedge} h_{p, w}^{\epsilon}$ and ${}^{\vee} h_{p, w}^{\epsilon}$ functions vary depending on the non-negativity and non-positivity, respectively, of their input vectors as follows:
    \scalebox{0.97}{$
    \begin{aligned}
        {}^{\wedge} h_{p, w}^{\epsilon}(y) &= 
        \begin{cases}
            (M_{0,w}^{\epsilon}(|y|_{+}^2))^{0.5} - \epsilon^{0.5} & \text{if $y_i \geq 0$, $\forall i \in \intv{1}{n}$ } \\
            \epsilon^{0.5} - (M_{p,w}^{\epsilon}(|y|_{-}^2))^{0.5} & \text{otherwise}
        \end{cases} \\
        {}^{\vee} h_{p, w}^{\epsilon}(y) &= 
        \begin{cases}
            \epsilon^{0.5} - (M_{0,w}^{\epsilon}(|y|_{-}^2))^{0.5} & \text{if $y_i \leq 0$, $\forall i \in \intv{1}{n}$ } \\
            (M_{p,w}^{\epsilon}(|y|_{+}^2))^{0.5} - \epsilon^{0.5} & \text{otherwise}.
        \end{cases}
    \end{aligned}
    $}
    The value of $\epsilon$, and the parameters associated with the generalized means, $p$, and $w$, can be chosen by considering the underlying problem as follows:
\end{rem}
\begin{itemize}
    \item 
    For any $\epsilon > 0$, ${}^{\wedge} h_{p, w}^{\epsilon}$ and ${}^{\vee} h_{p, w}^{\epsilon}$ functions demonstrate $\mathcal{C}^1$-smoothness,
    while preserving the soundness and completeness properties of the D-GMSR (see Lemma \ref{lem:sac}). It is desirable to select an arbitrarily small positive value for $\epsilon$ to ensure that the values of $M_{0,w}^{\epsilon}$ and $M_{p,w}^{\epsilon}$ functions closely approximate the weighted geometric mean and weighted power mean of their input.
    \item While $M_{p,w}^{0}(z)$ represents a weighted arithmetic mean for $p=1$,
    $(M_{p,w}^{\epsilon}(|z|_+^2))^{0.5} \to \max (z)$ and $(M_{p,w}^{\epsilon}(|z|_-^2))^{0.5} \to \min (z)$ as $p \to \infty$.
    Hence, users have the option to make the 
    values of the ${}^{\wedge} h_{p, w}^{\epsilon}$ and
    ${}^{\vee} h_{p, w}^{\epsilon}$
    functions as similar to the 
    values of the $\min$ and $\max$ functions under some conditions (see Lemma \ref{lem:limit}).
    \item The weight vector $w \in \mathbb{Z}_{++}^n$ allows users to express preferences and priorities, indicating the relative importance of elements in the vector $z \in \mathbb{R}^n_+$.
\end{itemize}

\begin{definition} \label{def:dgsar} (D-GMSR) D-GMSR is defined by replacing $\min$ and $\max$ functions of D-SR in Definition \ref{def:dtsr} with ${}^{\vee} h_{p, w}^{\epsilon}(y)$ and ${}^{\wedge} h_{p, w}^{\epsilon}(y)$ functions, respectively. The complete D-GMSR is defined as follows:
\scalebox{0.91}{
    $
    \begin{aligned}
        &\Gamma^{\mu} (x,k) := f (x_k)\\
        &\Gamma^{\neg \varphi}_{\bm{\epsilon}, \bm{p}, \bm{w}} (x,k) := - \Gamma_{\bm{\epsilon}, \bm{p}, \bm{w}}^\varphi (x,k)\\
        &\Gamma_{\bm{\epsilon}, \bm{p}, \bm{w}}^{\varphi_1 \wedge \varphi_2} (x,k) := {}^{\wedge} h_{p_1, w_1}^{\epsilon_1} ( (\Gamma^{\varphi_1}_{\bm{\epsilon}, \bm{p}, \bm{w}} (x,k), \Gamma^{\varphi_2}_{\bm{\epsilon}, \bm{p}, \bm{w}} (x,k)) )\\
        &\Gamma_{\bm{\epsilon}, \bm{p}, \bm{w}}^{\varphi_1 \vee \varphi_2} (x,k) := {}^{\vee} h_{p_1, w_1}^{\epsilon_1} ( (\Gamma^{\varphi_1}_{\bm{\epsilon}, \bm{p}, \bm{w}} (x,k), \Gamma^{\varphi_2}_{\bm{\epsilon}, \bm{p}, \bm{w}} (x,k)) )\\
        &\Gamma_{\bm{\epsilon}, \bm{p}, \bm{w}}^{\varphi_1 \implies \varphi_2} (x,k) := {}^{\vee} h_{p_1, w_1}^{\epsilon_1} ( (-\Gamma^{\varphi_1}_{\bm{\epsilon}, \bm{p}, \bm{w}} (x,k), \Gamma^{\varphi_2}_{\bm{\epsilon}, \bm{p}, \bm{w}} (x,k)) )\\
        &\Gamma_{\bm{\epsilon}, \bm{p}, \bm{w}}^{ \bm{F}_{\intv{a}{b}} \varphi}  (x,k) := {}^{\vee} h_{p_1, w_1}^{\epsilon_1}(( 
        \Gamma^{ \varphi }_{\bm{\epsilon}, \bm{p}, \bm{w}}(x, k+a), 
        \\&\hspace{3.958cm}
        \Gamma^{ \varphi }_{\bm{\epsilon}, \bm{p}, \bm{w}}(x, k+a+1), \\
        &\hspace{3.9cm}\dots, \Gamma^{ \varphi }_{\bm{\epsilon}, \bm{p}, \bm{w}}(x, k+b) ))\\
        &\Gamma_{\bm{\epsilon}, \bm{p}, \bm{w}}^{ \bm{G}_{\intv{a}{b}} \varphi}  (x,k) := {}^{\wedge} h_{p_1, w_1}^{\epsilon_1}(( 
        \Gamma^{ \varphi }_{\bm{\epsilon}, \bm{p}, \bm{w}}(x, k+a), 
        \\&\hspace{3.958cm}
        \Gamma^{ \varphi }_{\bm{\epsilon}, \bm{p}, \bm{w}}(x, k+a+1),\\
        &\hspace{3.9cm}\dots, \Gamma^{ \varphi }_{\bm{\epsilon}, \bm{p}, \bm{w}}(x, k+b) ))\\
        &\Gamma_{\bm{\epsilon}, \bm{p}, \bm{w}}^{ \varphi_1 \bm{U}_{\intv{a}{b}} \varphi_2}  (x,k) := {}^{\vee} h_{p_1, w_1}^{\epsilon_1}(( {}^{\wedge} h_{p_2, w_2}^{\epsilon_2}(( \Gamma_{\bm{\epsilon}, \bm{p}, \bm{w}}^{ \bm{G}_{\intv{a}{a}} \; \varphi_1}(x, k), \\
        &\hspace{5.915cm}\Gamma_{\bm{\epsilon}, \bm{p}, \bm{w}}^{ \varphi_2 }(x, k+a) )), \\
        &\hspace{4.46cm}{}^{\wedge} h_{p_2, w_2}^{\epsilon_2}(( \Gamma_{\bm{\epsilon}, \bm{p}, \bm{w}}^{ \bm{G}_{\intv{a}{a+1}} \; \varphi_1}(x, k), \\
        &\hspace{5.915cm}\Gamma_{\bm{\epsilon}, \bm{p}, \bm{w}}^{ \varphi_2 }(x, k+a+1) )), \\
        &\hspace{4.75cm}\vdots \\
        &\hspace{4.46cm}{}^{\wedge} h_{p_2, w_2}^{\epsilon_2}(( \Gamma_{\bm{\epsilon}, \bm{p}, \bm{w}}^{ \bm{G}_{\intv{a}{b}} \; \varphi_1}(x, k), \\
        &\hspace{5.92cm}\Gamma_{\bm{\epsilon}, \bm{p}, \bm{w}}^{ \varphi_2 }(x, k+b) ))))
    \end{aligned}
    $}
where $\bm{\epsilon}$, $\bm{p}$, and $\bm{w}$ represent sets consisting of an arbitrary number of $\epsilon_i \in \mathbb{R}_{++}$, $p_i \in \mathbb{Z}_{++}$, and $w_i \in \mathbb{Z}_{++}^n$, respectively, required for the construction of ${}^{\wedge} h_{p_i, w_i}^{\epsilon_i}$ and ${}^{\vee} h_{p_i, w_i}^{\epsilon_i}$ functions, which in turn construct the corresponding specification.
For instance, let $\varphi := \bm{F}_{\intv{1}{5}} ((f_1(x) \geq 0) \wedge (f_2(x) \geq 0))$. Then, 
the cardinality of the sets $\bm{\epsilon}$, $\bm{p}$ and $\bm{w}$ used in $\Gamma_{\bm{\epsilon}, \bm{p}, \bm{w}}^\varphi (x, 0)$ is $2$, where $w_1 \in \mathbb{Z}_{++}$ and $w_2 \in \mathbb{Z}_{++}^5$.
\end{definition}

\subsection{Properties of the D-GMSR}
\vspace{-0.25cm}
In Lemma \ref{lem:smth} and \ref{lem:sac}, we assume $\epsilon_i \in \mathbb{R}_{++}$, $p_i \in \mathbb{Z}_{++}$, $w_i \in \mathbb{Z}_{++}^{n_i}$ for all $i \in \mathbb{Z}_{++}$ and do not restate this for brevity.
\vspace{-0.25cm}
\begin{lemma} \label{lem:smth} ($\mathcal{C}^1$-smoothness)
    $\Gamma^\varphi_{\bm{\epsilon}, \bm{p}, \bm{w}} ( \cdot, k) \in \mathcal{C}^1(\mathbb{R}^n)$ 
    for any  finite time STL specification $\varphi$ and time step $k \in \mathbb{Z}_+$.
\end{lemma}
\vspace{-0.25cm}
\begin{pf}
    Given that Assumption \ref{asm:pred_smth} requires all predicate functions to be $\mathcal{C}^1$-smooth,
    finite time STL specification $\varphi$ and time step $k \in \mathbb{Z}_+$ if and only if each
    ${}^{\wedge} h_{p_i, w_i}^{\epsilon_i}$ and ${}^{\vee} h_{p_i, w_i}^{\epsilon_i}$ functions that construct the STL specification $\varphi$ are $\mathcal{C}^1$-smooth on $\mathbb{R}^n$.
    $|y|_{\Hsquare}^2 \in \mathbb{R}_+^n$ implies that $M_{p_i, w_i}^{\epsilon_i}(|y|_{\Hsquare}^2), M_{0, w_i}^{\epsilon_i}(|y|_{\Hsquare}^2) \in [\epsilon, \infty ]$. 
    \fbs{Note that $\nabla_x | x |_{\Hsquare}^2 = 2 | x |_{\Hsquare} $ for $x \in \mathbb{R}^n$, which implies $| \cdot |_{\Hsquare}^2 \in \mathcal{C}^1  (\mathbb{R}^n)$
    Then}
    $M_{p_i, w_i}^{\epsilon_i}, M_{0, w_i}^{\epsilon_i} \in \mathcal{C}^{\infty} (\mathbb{R}_+^n)$ implies that $M_{p_i, w_i}^{\epsilon_i}( |\cdot|_{\Hsquare}^2 ), M_{0, w_i}^{\epsilon_i}( |\cdot|_{\Hsquare}^2 ) \in \mathcal{C}^1(\mathbb{R}^n)$.
    Therefore, ${}^{\wedge} h_{p_i, w_i}^{\epsilon_i}, {}^{\vee} h_{p_i, w_i}^{\epsilon_i}$ $\in \mathcal{C}^1 (\mathbb{R}^n)$.
\end{pf}
\vspace{-0.25cm}
\begin{lemma} \label{lem:sac} (Soundness and Completeness) 
    For any discrete-time signal $x$, time step $k \in \mathbb{Z}_+$, and finite time STL specification $\varphi$, $(x, k) \models \varphi \iff \Gamma^\varphi_{\bm{\epsilon}, \bm{p}, \bm{w}} (x, k) \geq 0$.
\end{lemma}
\begin{pf}
    Given that D-SR in Definition \ref{def:dtsr} is both sound and complete \cite{gilpin2020smooth}, and ${}^{\vee} h_{p_i, w_i}^{\epsilon_i}(y) = -{}^{\wedge} h_{p_i, w_i}^{\epsilon_i}(-y)$, \bYY{analogous to $\max(y) = -\min(-y)$, proofs of the following assertions establish the soundness and completeness of D-GMSR.} 
    \begin{itemize}
        \item ${}^{\wedge} h_{p_i, w_i}^{\epsilon_i}(y) < 0\iff\min(y) < 0$,
        \item ${}^{\wedge} h_{p_i, w_i}^{\epsilon_i}(y) > 0\iff\min(y) > 0$.
    \end{itemize}
    One can verify from (\ref{eq:key_fcn}) that 
    \begin{itemize}
        \item $M_{0, w_i}^{\epsilon_i}(|y|_{+}^2) > \epsilon \iff y_j > 0$ for all $j \in \intv{1}{n}$,
        \item $M_{p_i, w_i}^{\epsilon_i}(|y|_{-}^2) > \epsilon \iff \exists j \in \intv{1}{n}$ such that  $y_j < 0$,
        \item $M_{0, w_i}^{\epsilon_i}(|y|_{+}^2) > \epsilon \implies M_{p_i, w_i}^{\epsilon_i}(|y|_{-}^2)  = \epsilon$,
        \item $M_{p_i, w_i}^{\epsilon_i}(|y|_{-}^2) > \epsilon \implies M_{0, w_i}^{\epsilon_i}(|y|_{+}^2) = \epsilon$.
    \end{itemize}
    Therefore,
    \begin{align*}
        \min(y) < 0 &\iff \exists j  \in \intv{1}{n} \text{ such that } y_j < 0 \\
        &\iff M_{p_i, w_i}^{\epsilon_i}(|y|_{-}^2) > \epsilon \\
        &\iff M_{0, w_i}^{\epsilon_i}(|y|_{+}^2) - M_{p_i, w_i}^{\epsilon_i}(|y|_{-}^2) < 0 \\
        &\iff {}^{\wedge} h_{p_i, w_i}^{\epsilon_i}(y) < 0
    \end{align*}
    \begin{align*}
        \min(y) > 0 &\iff y_j > 0 \text{ for all } j  \in \intv{1}{n} \\
        &\iff M_{0, w_i}^{\epsilon_i}(|y|_{+}^2) > \epsilon \\
        &\iff M_{0, w_i}^{\epsilon_i}(|y|_{+}^2) - M_{p_i, w_i}^{\epsilon_i}(|y|_{-}^2) > 0 \\
        &\iff {}^{\wedge} h_{p_i, w_i}^{\epsilon_i}(y) > 0
    \end{align*}
\end{pf}
\begin{lemma} \label{lem:limit} (Limit behavior)  
    \begin{itemize}
        \item $\lim_{p \to \infty} \; {}^{\wedge} h_{p, w}^{\epsilon}(y) = \min (y) + \epsilon^{0.5}$ if $\min (y) < 0$ and $\epsilon^{0.5} < | \min (y) |$.
        \item  $\lim_{p \to \infty} \; {}^{\vee} h_{p, w}^{\epsilon}(y) = \max (y) - \epsilon^{0.5}$ if $\max (y) > 0$ and $\epsilon^{0.5} < \max (y)$.
    \end{itemize}
\end{lemma}
\begin{pf}
    If $z \in \mathbb{R}_+$, then $M_{\infty, w}^{0}(z) = \max (z)$ for any $w \in \mathbb{R}_{++}^n$ \cite{shniad1948convexity}. If $\epsilon^{0.5} < \max (z)$, then $M_{\infty, w}^{\epsilon}(z) = \max (z)$. 
    Therefore, if $\min (y) < 0$ and $\epsilon^{0.5} < |\min (y)|$, then $\lim_{p \to \infty} \; {}^{\wedge} h_{p, w}^{\epsilon}(y) = 
    \epsilon^{0.5} - \lim_{p \to \infty}  \Big( M_{p, w}^{\epsilon}( |y|_{-}^2 ) \Big)^{1/2} = \epsilon^{0.5} - ( \max  (|y|_{-}^2) )^{1/2} = \epsilon^{0.5} - | \max (|y_i|_{-}) | =  \min (y) + \epsilon^{0.5}$.
    
    Since ${}^{\vee} h_{p, w}^{\epsilon}(y) := -{}^{\wedge} h_{p, w}^{\epsilon}(-y)$, $\lim_{p \to \infty} \; {}^{\vee} h_{p, w}^{\epsilon}(y) =  -\min (-y) - \epsilon^{0.5} = \max (y) - \epsilon^{0.5} $ if $\max (y) > 0$ and $\epsilon^{0.5} < \max (y)$.
\end{pf}

\subsection{Gradient properties of the D-GMSR} \label{sec:grad_dgsar}
As discussed in Section \ref{sec:grad_dssr}, the robustness value of any STL specification is determined by the value of a specific predicate function at a single time instant in D-SR since all predicate functions included in the STL specification and their values at all times are composed with the $\min$ and $\max$ functions.
Accordingly, when maximizing the robustness value of a specification using D-SSR and a gradient-based optimization method, typically only the gradient of a specific predicate function at a single time instant is considered notably in each iteration.
Conversely, key functions of D-GMSR, denoted as ${}^{\wedge} h_{p, w}^{\epsilon}$ and ${}^{\vee} h_{p, w}^{\epsilon}$, are formed through a specific composition of generalized means in (\ref{eq:dgsar_key}). This construction enables the simultaneous consideration of values from multiple predicate functions across multiple time instances, addressing the issues of locality and masking observed in D-SR.
Accordingly, the gradient of an STL specification is the weighted summation of the gradient of multiple predicate functions at multiple time instances.
These characteristics of key functions of D-GMSR lead to the alleviation of certain undesirable optimization results observed in D-SSR, where the STL specification is not met, as demonstrated by an example in Section \ref{sec:num_sim}.

To demonstrate the characteristic of key functions of D-GMSR, suppose that $f_i: \mathbb{R}^n \to \mathbb{R}$ is a function of $x \in \mathbb{R}^n$ for all $i \in \intv{1}{m}$ and $y = (f_1(x), f_2(x), \dots, f_m(x) ) $, then 
\begin{align*}
    \nabla_x {}^{\wedge} h_{p, w}^{\epsilon}(&y) = c_0^+ \nabla_x M_{0, w}^{\epsilon}(|y|_{+}^2) - c_p^- \nabla_x M_{p, w}^{\epsilon}(|y|_{-}^2), \\
    \nabla_x {}^{\vee} h_{p, w}^{\epsilon}(&y) = c_p^+ \nabla_x M_{p, w}^{\epsilon}(|y|_{+}^2) - c_0^- \nabla_x M_{0, w}^{\epsilon}(|y|_{-}^2),
\end{align*}
where
\begin{align*}
    \nabla_x M_{p, w}^{\epsilon}(|y|_{\square}^2) &= c_{p,m}^{\square} \sum_{i=1}^m w_i |f_i(x)|_{\square}^{2 p - 1} \nabla_x f_i(x), \\
    \nabla_x M_{0, w}^{\epsilon}(|y|_{\square}^2) &= c_{0,m}^{\square} \sum_{i=1}^m   w_i  \frac{1}{|f_i(x)|_{\square}} \nabla_x f_i(x),
\end{align*}
\begin{align*}
    c_0^{\square} &= \frac{1}{2} \Big( M_{0, w}^{\epsilon}(|y|_{\square}^2) \Big)^{-1/2},\\
    c_p^{\square} &= \frac{1}{2} \Big( M_{p, w}^{\epsilon}(|y|_{{\square}}^2) \Big)^{-1/2},\\
    c_{0,m}^{\square} &= \frac{2 \fbs{\prod_{j = 1}^m} |f_j(x)|_{\square}^{2w_i} }{
    \boldsymbol{1}^T w \Big( M_{0, w}^{\epsilon}(|y|_{\square}^2) \Big)^{\boldsymbol{1}^T w-1} 
    }, \\
    c_{p,m}^{\square} &=  \frac{2p}{
    p \boldsymbol{1}^T w \Big( M_{p, w}^{\epsilon}(|y|_{\square}^2) \Big)^{p-1} 
    }.
\end{align*}
Consequently,
\begin{align*}
    \nabla_x  {}^{\wedge} h_{p, w}^{\epsilon}( y ) &= \fbs{\nabla_x  {}^{\wedge} h_{p, w}^{\epsilon}( (f_1(x), f_2(x), \dots, f_m(x) ) )} \\
    &= \sum_{i=1}^m w_i a_i \nabla_x f_i(x) \text{ and}\\
    \nabla_x  {}^{\vee} h_{p, w}^{\epsilon}( y ) &= \fbs{\nabla_x  {}^{\vee} h_{p, w}^{\epsilon}( (f_1(x), f_2(x), \dots, f_m(x) ) )} \\
    &= \sum_{i=1}^m w_i b_i \nabla_x f_i(x)    
\end{align*}
where
\begin{align*}
    a_i &= 
    \begin{cases}
        c_0^+ c_{0,m}^+ \displaystyle \frac{1}{ | f_i(x) |_+ } & \text{if } {}^{\wedge} h_{p, w}^{\epsilon}( y ) > 0\\[0.3cm]
        c_p^- c_{p,m}^-  -| f_i(x) |^{2p-1}_- & \text{otherwise}
    \end{cases}, \\
    b_i &= 
    \begin{cases}
        c_0^- c_{0,m}^- \displaystyle
        \frac{1}{ -| f_i(x) |_- } & \text{if } {}^{\vee} h_{p, w}^{\epsilon}( y ) < 0\\[0.3cm]
        c_p^+ c_{p,m}^+  | f_i(x) |^{2p-1}_+ & \text{otherwise}
    \end{cases}.
\end{align*}
Therefore,
\begin{itemize}
    \item If ${}^{\wedge} h_{p, w}^{\epsilon}( y ) > 0$, then 
    $\frac{a_i}{a_j} \propto \frac{ f_j(x)}{ f_i(x)}$ for all $i,j \in \intv{1}{m} $.

    \item If ${}^{\wedge} h_{p, w}^{\epsilon}( y ) < 0$, then 
    $ \frac{a_i}{a_j} \propto \frac{ |f_i(x)|^{2p-1} }{ |f_j(x)|^{2p-1} }$ for all $i,j \in \{ k : f_k < 0, k \in \intv{1}{m} \}$, and $a_i = 0$ for all $i \in \{ k : f_k \geq 0, k \in \intv{1}{m} \}$.

    \item If ${}^{\vee} h_{p, w}^{\epsilon}( y ) < 0$, then 
    $\frac{b_i}{b_j} \propto \frac{ |f_j(x)| }{ |f_i(x)|}$ for all $i,j \in \intv{1}{m}$.

    \item If ${}^{\vee} h_{p, w}^{\epsilon}( y ) > 0$, then 
    $\frac{b_i}{b_j} \propto \frac{ f_i(x)^{2p-1}}{ f_j(x)^{2p-1}}$ for all $i,j \in \{ k : f_k > 0, k \in \intv{1}{m} \}$, and $b_i = 0$ for all $i \in \{ k : f_k \leq 0, k \in \intv{1}{m} \}$.
\end{itemize}
Recall that for the $\min$ and ${}^{\wedge} h_{p, w}^{\epsilon}$ functions to be positive, all their variables must be positive.
The derivative of the $\fbs{\widetilde{\min}_{\kappa}}$ function yields almost zero for all the variables that do not have the minimum value. On the other hand, 
if the values of all the variables are negative, then 
the derivative of the ${}^{\wedge} h_{p, w}^{\epsilon}$ function with respect to each variable is non-zero and it increases as the value of the variable decreases.
In other words, the variable that is farthest from being positive receives greater consideration during the optimization process while other negative-valued variables are still not ignored.
If there exist both positive and negative-valued variables, then the derivative with respect to negative-valued variables increases as the value of the variable decreases 
and the derivative of the function with respect to non-negative-valued variables becomes zero.
In other words, only the variables causing the negative yield of the ${}^{\wedge} h_{p, w}^{\epsilon}$ function are taken into account.

For the $\max$ and ${}^{\vee} h_{p, w}^{\epsilon}$ functions to be positive, only one of their variables must be positive.
The derivative of the $\fbs{\widetilde{\max}_{\kappa}}$ function yields almost zero for all the variables that do not have the maximum value. On the other hand,  
if the values of all the variables are negative, then 
the derivative of the ${}^{\vee} h_{p, w}^{\epsilon}$ function with respect to each variable is non-zero and it increases as the value of the variable increases.
In other words, the variable that is closest to being positive receives greater consideration during the optimization process while other negative-valued variables are still not ignored.
If there exist both positive and negative-valued variables, then the derivative with respect to positive-valued variables increases as the value of the variable increases 
and the derivative of the function with respect to non-positive-valued variables becomes zero.
In other words, only the variables causing the positive yield of the ${}^{\wedge} h_{p, w}^{\epsilon}$ function are taken into account.

As a result, unlike the $\fbs{\widetilde{\min}_{\kappa}}$ and $\fbs{\widetilde{\max}_{\kappa}}$ functions, 
while the values of the ${}^{\wedge} h_{p, w}^{\epsilon}$ and ${}^{\vee} h_{p, w}^{\epsilon}$ functions are negative, their derivatives with respect to the variables whose variations can make these function positive are never zero.
Hence, certain 
undesirable optimization results observed in D-SSR, where the STL specification is not met, are mitigated, as demonstrated by an example in Section \ref{sec:num_sim}.
\begin{example}
Similar to the example presented in Figure \ref{fig:dssr}, we illustrate the variations in the predicate function values during the maximization of $\Gamma^{\varphi_1 \wedge \varphi_2 \wedge \dots \wedge \varphi_5}_{\bm{\epsilon}, \bm{p}, \bm{w}} (x,k)$ and $\Gamma^{\varphi_1 \vee \varphi_2 \vee \dots \vee \varphi_5}_{\bm{\epsilon}, \bm{p}, \bm{w}} (x,k)$ functions in Figure \ref{fig:dgsar}
to provide numerical evidence for the statements regarding the derivative of ${}^{\wedge} h_{p, w}^{\epsilon}$ and $ {}^{\vee} h_{p, w}^{\epsilon}$ functions for $\epsilon = 10^{-8}$, $p = 1$ and $w = \bm{1}$.
Here, $\Gamma^{\varphi_i}(x,k) = f_i(x_k) \in \mathbb{R}$ and $f_i : \mathbb{R} \to \mathbb{R}$ for all $i \in \intv{1}{5}$. Therefore, 
$\Gamma^{\varphi_1 \wedge \varphi_2 \wedge \dots \wedge \varphi_5}_{\bm{\epsilon}, \bm{p}, \bm{w}} (x,k) = {}^{\wedge} h_{p, w}^{\epsilon} ( y ) $ and $\Gamma^{\varphi_1 \vee \varphi_2 \vee \dots \vee \varphi_5}_{\bm{\epsilon}, \bm{p}, \bm{w}} (x,k) = {}^{\vee} h_{p, w}^{\epsilon} ( y ) $, where $y=( f_1(x_k), f_2(x_k), \dots, f_5(x_k) )$.
Alternatively, maximizing $\Gamma^{ \bm{G}_{[1 : 5]} \varphi}_{\bm{\epsilon}, \bm{p}, \bm{w}} (x,0)$ and $\Gamma^{ \bm{F}_{[1 : 5]} \varphi}_{\bm{\epsilon}, \bm{p}, \bm{w}} (x,0)$ functions, respectively, where $\Gamma^{\varphi}(x,k) = f(x_k) \in \mathbb{R}$ for all $k \in \intv{1}{5}$ and $f : \mathbb{R} \to \mathbb{R}$, leads to similar changes in predicate function values across different time instances because $\Gamma^{\bm{G}_{[1 : 5]} \varphi}_{\bm{\epsilon}, \bm{p}, \bm{w}} (x,0) = {}^{\wedge} h_{p, w}^{\epsilon} ( y )$ and $\Gamma^{\bm{F}_{[1 : 5]} \varphi}_{\bm{\epsilon}, \bm{p}, \bm{w}} (x,0) = {}^{\vee} h_{p, w}^{\epsilon} ( y ) $, where $y=( f(x_1), f(x_2), \dots, f(x_5)  )$.

For ${}^{\wedge} h_{p, w}^{\epsilon}$ function, 
when there exist both positive and negative-valued variables, the derivative of the function with respect to its variables causes only the value of the negative-valued variables to increase such that their value will become equal to each other at $0$, regardless of their initial values.
When there exists only positive-valued variables, all variables increase such that their values will become equal to each other since the variable that has a smaller value increases more.
For ${}^{\vee} h_{p, w}^{\epsilon}$ function, when there exist only negative-valued variables, all of them increase such that the variable that is closest to being positive increases more. When there exist both positive and negative-valued variables, only positive-valued variables increase such that the variable with higher value increases more than other positive-valued variables.
These results validate the statements regarding the derivatives of ${}^{\wedge} h_{p, w}^{\epsilon}$ and ${}^{\vee} h_{p, w}^{\epsilon}$ functions. 
\begin{figure}
\centerline{\includegraphics[scale=0.43]{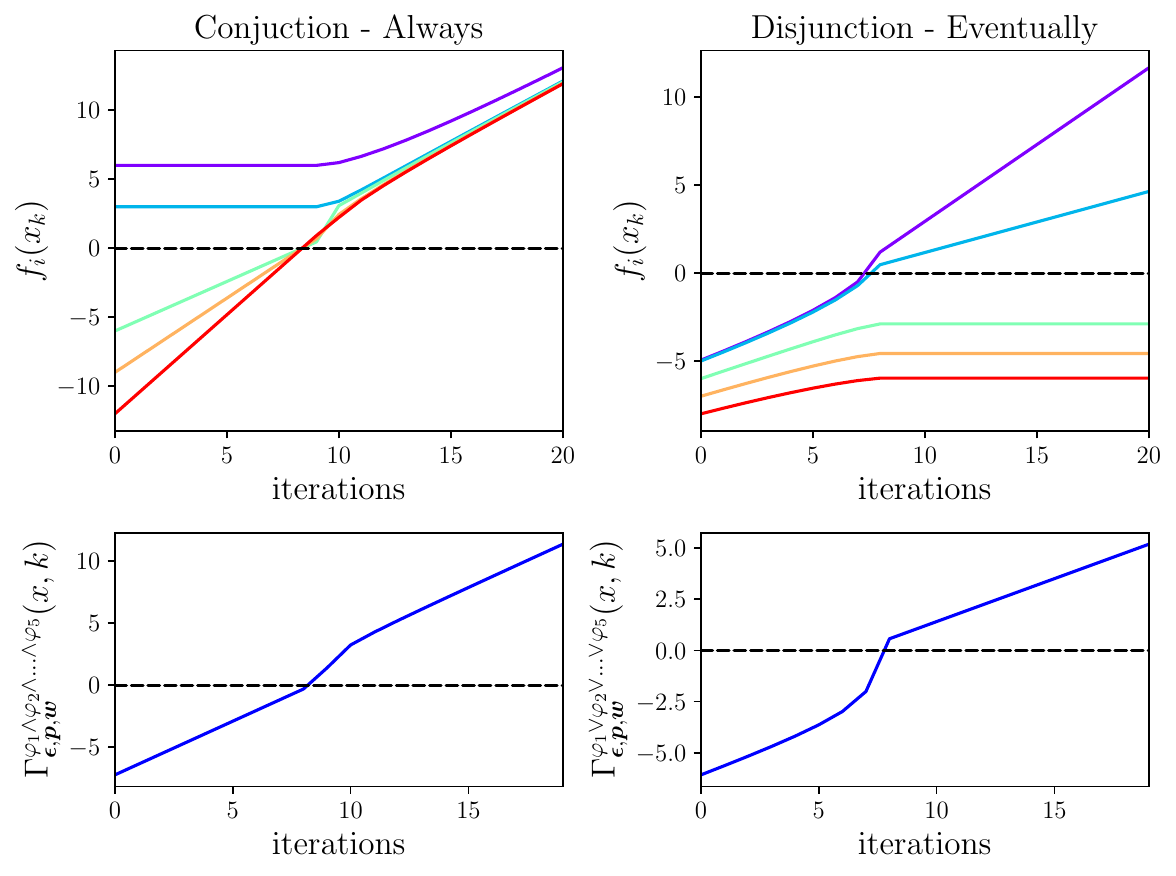}}
\caption{
The changes in the values of $f_i(x_k)$ (or $f(x_i)$) and the corresponding robustness values while maximizing the \textit{conjunction} (or \textit{always}) and \textit{disjunction} (or \textit{eventually}) operators using D-GMSR are presented on the left and right sides of the figure, respectively.
\bYY{The color gradient, ranging from purple to red in the figures above, represents the values of the functions from $f_1(x_k)$ to $f_5(x_k)$ (or $f(x_1)$ to $f(x_5)$), respectively.}
}
\label{fig:dgsar}
\end{figure}
\end{example}
\section{Numerical Simulations} \label{sec:num_sim}
Sequential convex programming (SCP) is a prominent direct method to solve trajectory optimization problems \cite{mao2016successive, malyuta2021advances, malyuta2022convex, ctcs2024}. 
SCP algorithms are also utilized to initialize indirect methods for trajectory optimization \cite{tang2018fuel, spada2023direct}.
The real-time capability of this approach is demonstrated for powered descent guidance problems in \cite{Elango2022, Kamath2023}. 
SCP algorithms also demonstrate good performance in solving optimal control problems with state-triggered constraints, which are a subset of STL constraints \cite{szmuk2019real,szmuk2020successive,reynolds2020dual}.

For solving the trajectory optimization problems we present, we use the successive convexification framework \cite{ctcs2024}, \fbs{which} combines multiple-shooting for the exact discretization of the continuous-time dynamics \cite{bock1984multiple, quirynen2015multiple}, $\ell_1$ exact penalization for the nonconvex constraints \cite[Chap. 17]{nocedal2006numerical} and a convergence-guaranteed SCP algorithm, prox-linear method \cite{drusvyatskiy2019efficiency}, for solving the resulting optimization problem. 
For modeling and solving each convex subproblem encountered in every iteration of the prox-linear algorithm, we use CVXPY \cite{diamond2016cvxpy} with either ECOS \cite{domahidi2013ecos} or MOSEK \cite{aps2019mosek}.
Customizable conic solvers such as the proportional-integral projected gradient (PIPG) algorithm \cite{yu2022proportional, yu2022extrapolated} could also be employed for efficient implementations.
We employ the adaptive weight update algorithm \cite[Alg. 2.2]{cartis2011evaluation} to dynamically determine the proximal term weight within the prox-linear algorithm. 
The convergence guarantee of the prox-linear method, achieved through an appropriate selection of the proximal term weight, ensures that the converged solutions, which are feasible with respect to the finite-dimensional nonconvex optimization problem, are also Karush-Kuhn-Tucker (KKT) points \cite[Thm. 19]{ctcs2024}.
The numerical performance of the successive convexification framework is demonstrated for applications ranging from GPU-accelerated trajectory optimization for 6-DoF rocket landing \cite{chari2024fast} and nonlinear model predictive control (NMPC) for obstacle avoidance \cite{nmpc2024}, to trajectory optimization for 6-DoF aircraft approach and landing \cite{aircraft2024}.
\begin{rem} \label{rem:iscv}
    A reformulation technique is proposed for the successive convexification framework in \cite{ctcs2024} to ensure the satisfaction of the path constraints in continuous time; however, \fbs{it covers only \textit{always} specification in STL. The application of this approach to STL specifications in general has not been explored.} Therefore, given that we solely considered STL specifications at node points, we impose path constraints exclusively at node points for consistency.
\end{rem}
\vspace{-0.15cm}
\subsection{Locality \& Masking}
\begin{figure}
\centerline{\includegraphics[scale=0.43]{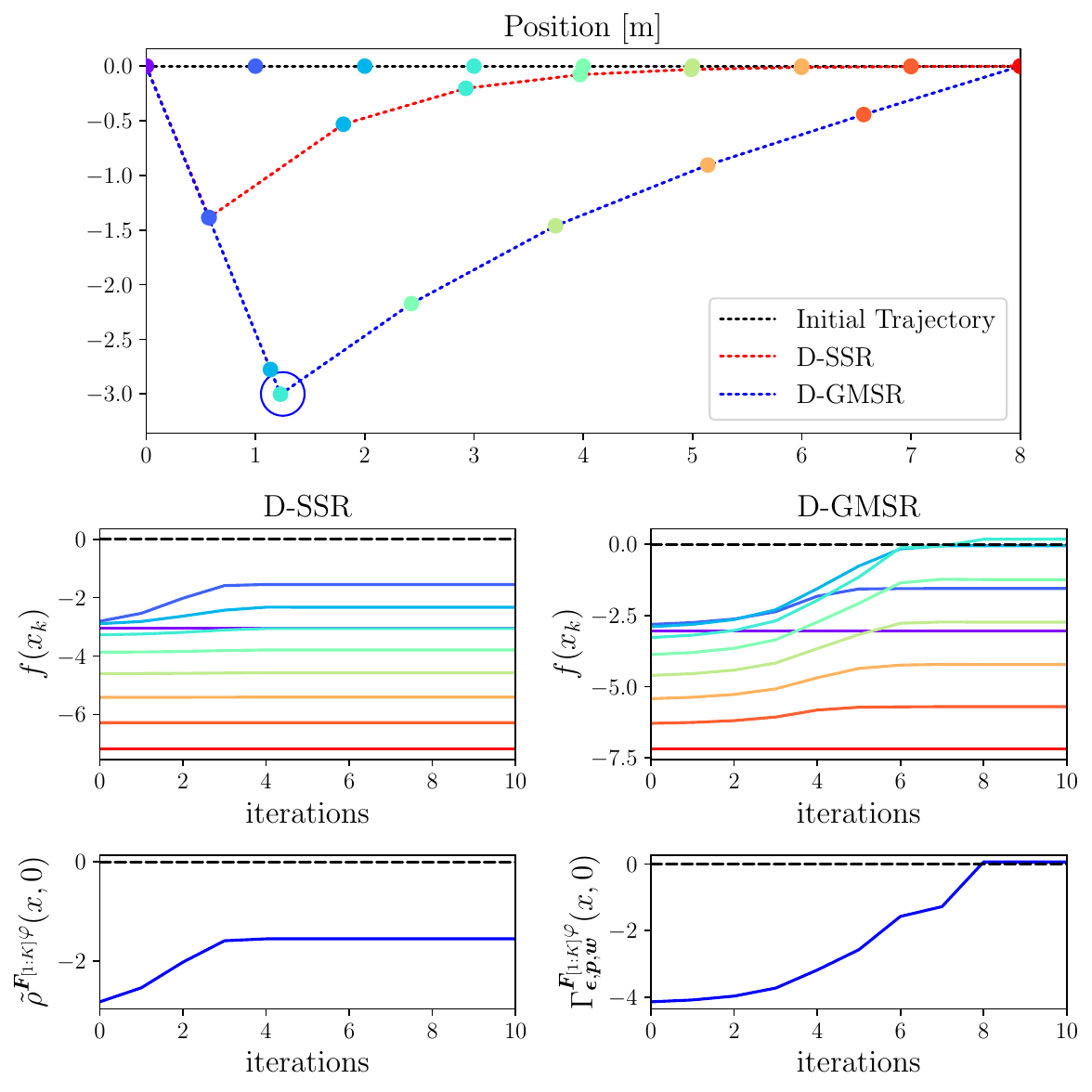}}
\caption{ 
The initial and converged trajectories obtained using D-SSR and D-GMSR, along with the changes in predicate function for each node point, and the corresponding robustness values throughout the iterations, are illustrated.
In the position plot, the blue circle shows the area that needs to be visited. 
The black dotted curve represents the initial trajectory, while the red and blue dotted curves show the converged trajectories obtained using D-SSR and D-GMSR, respectively. 
Node points of each trajectory, from $x_1$ to $x_9$, are represented with varying colors ranging from purple to red.
Node points are connected by dotted curves for clarity.
Similarly, in the central figures, the color gradient, ranging from purple to red, represents the value of the predicate function for each node point from $f(x_1)$ to $f(x_9)$, respectively.
}
\label{fig:disj_or_cvx_agm_all}
\end{figure}

\vspace{-0.15cm}
We first solve a simple trajectory optimization problem in (\ref{eq:ct_opt_cont_3}) to illustrate an 
undesirable optimization result, where the STL specification is not met due to the properties of $\fbs{\widetilde{\min}_{\kappa}}$ and $\fbs{\widetilde{\max}_{\kappa}}$ functions in D-SSR, followed by a demonstration of D-GMSR's performance on the same problem.
The position and control input of the vehicle at time $k$ are represented as $x_k$ and $u_k$, respectively. The control input of the vehicle is upper-bounded and there exist boundary conditions on the initial and final positions of the vehicle.
The objective is for the vehicle to \textit{eventually} reach a circular area centered at $p_w$ with a radius of $r_w$. The corresponding predicate function and the specification are defined as $f(x) := r_w - \| x - p_w \|$ and $\varphi := (f(x) \geq 0)$, respectively.
The solutions of the optimization problems utilize D-SSR for $\kappa = 25$ \cite{gilpin2020smooth} and D-GMSR for $\epsilon=10^{-8}$, $p=1$, $w=\bm{1}$ 
are presented in Fig. (\ref{fig:disj_or_cvx_agm_all}).
\begin{empheq}[box=\fbox]{equation}
    \begin{aligned} \label{eq:ct_opt_cont_3}
        \underset{x_k, u_k}{\text{maximize}} \quad & \tilde{\rho}^{ \bm{F}_{\intv{1}{K}} \varphi}  (x,0) \text{ or } \Gamma^{ \bm{F}_{\intv{1}{K}} \varphi}_{\bm{\epsilon}, \bm{p}, \bm{w}}  (x,0) \\
        \text{subject to} \quad 
        &  x_{k+1} = x_k + u_k, \; \forall k \in \intv{1}{K-1}  \\
        &  x_1 = (0, 0)  \\
        &  x_K = (8, 0)  \\
        & \| u_k \| \leq 1.5, \; \forall k \in \intv{1}{K}
    \end{aligned}
\end{empheq}
The total number of nodes ($K$) is set as $9$ and the initial trajectory for the optimization problem is set as the black dotted curve in Fig. \ref{fig:disj_or_cvx_agm_all}. In order to satisfy the STL specification, at least one of the node points ($x_k$) must be inside of the blue circle in Fig. \ref{fig:disj_or_cvx_agm_all}. 
For the initial trajectory, the second node point is the closest node point to the center of the blue circle, thus $\partial_{x_k} \tilde{\rho}^{ \bm{F}_{\intv{1}{K}} \varphi}  (x,0) \approx 0 $ for all $k \in \intv{1}{K} \setminus \{ 2 \} $. 
In other words, when the D-SSR is used, the optimizer only forces the second node point to converge to the center of the blue circle. 
However, due to the control-upper bound, the distance between two node points cannot be greater than $1.5$. 
Therefore, $x_2$ cannot converge to the inside of the blue circle since $x_1 = (0,0)$. 
As a result, the initial trajectory converges to the trajectory presented via the red dotted curve in Fig. \ref{fig:disj_or_cvx_agm_all}.
Even though alternative node points could satisfy the specification, D-SSR's focus solely on the second node point prevents the satisfaction of the specification.

In D-GMSR, however, for the initial trajectory, we have
\begin{equation*}
    \partial_{x_k} \Gamma^{ \bm{F}_{\intv{1}{K}} \varphi}_{\bm{\epsilon}, \bm{p}, \bm{w}}  (x,0) = a_k \partial_{x_k} f(x)
\end{equation*}
where $a_2 > a_k \gg 0$ for all $k \in \intv{1}{K} \setminus \{ 2 \} $. 
In other words, when D-GMSR is used, the optimizer prioritizes the convergence of the second node point towards the center of the blue circle but it also enforces other node points inversely proportional to their distances from the circular area. 
Therefore, despite the control upper bound preventing the convergence of the second node point within the blue circle, the optimization process ensures that other node points continue to converge towards the circle's center, eventually resulting in one of them reaching the center.
As a result, when D-GMSR is used, the initial trajectory converges to the trajectory presented via the blue dotted curve in Fig. \ref{fig:disj_or_cvx_agm_all}.

\vspace{-0.15cm}
\subsection{Quadrotor flight} \label{sec:qf}
\vspace{-0.15cm}
We address an optimal control problem for a point mass model of a quadrotor, incorporating STL specifications. In this scenario, the objective for the quadrotor is to fly to a designated destination while satisfying the constraints on speed, tilt angle, and thrust force. However, the vehicle's speed must not exceed a specified threshold \textit{until} it stops at the battery charging station for a consecutive number of time steps.

\subsubsection{Vehicle dynamics} 
The state and control input of the vehicle are defined as $x := (r, v)$ and $u := T$, where $r := (r^x, r^y, r^z) \in \mathbb{R}^3$, $v := (v^x, v^y, v^z) \in \mathbb{R}^3$, and $T := (T^x, T^y, T^z) \in \mathbb{R}^3$ are the position, velocity, and thrust force, respectively. The dynamic equation of the vehicle is defined as
\begin{align} \label{eq:quad_dyn_eq}
    \Dot{x}(t) &=  f(x(t), u(t))\\ 
    &= 
    \begin{bmatrix}
        v(t) \\
        T(t)/m - c_{\mathrm{d}} \| v(t) \| v(t) + g
    \end{bmatrix} \nonumber
\end{align}
where $m \in \mathbb{R}_{++}$ and $c_{\mathrm{d}} \in \mathbb{R}_{++}$ are the mass and the drag coefficient of the vehicle, respectively, and $g \in \mathbb{R}^3$ is the gravity vector.

We discretize the continuous-time dynamics using the multiple-shooting approach. The time is discretized over a finite grid as $0 = t_1 < t_2 < \dots < t_K = t_f$, $x_k$ and $u_k$ represent the state and control input at the node point $t_k$ for $k \in \intv{1}{K}$. The control input $u$ is parameterized via first-order hold and the continuous-time dynamics \eqref{eq:quad_dyn_eq} is discretized using multiple-shooting as follows:
\begin{align} \label{eq:quad_dyn_dt}
    x_{k+1} = F_k(x_k, u_k, u_{k+1})
\end{align}
for $k \in \intv{1}{K-1}$. See \cite{ctcs2024, nmpc2024} for further details.

The boundary conditions are specified as follows:
\begin{subequations} \label{eq:quad_bc}
\begin{gather}
    r_1 = r_{\mathrm{i}}, \; r_K = r_{\mathrm{f}}, \; v_1 = v_{\mathrm{i}}, \; v_K = v_{\mathrm{f}} \tag{\theequation a-d} \\
    \!\!\!\!\!\!\!\!\!\!\!\!\!\! u_1 = mg, \; u_K = mg \tag{\theequation e-f}
\end{gather}
\end{subequations}
The speed, the tilt angle, and the thrust force of the vehicle are constrained as follows:
\begin{subequations} \label{eq:quad_cons}
\begin{align}
    \| v_k \|_2 &\leq v_{\mathrm{max}} \;\;\; \forall k \in \intv{1}{K} \\
    \| u_k \|_2 \cos(\theta) &\leq e_1^{\top} u_k \;\;\; \forall k \in \intv{1}{K} \\
    \| T_k \|_2 &\leq T_{\mathrm{max}} \;\;\; \forall k \in \intv{1}{K}
\end{align}
\end{subequations}
\subsubsection{STL specifications} 
In this scenario, the vehicle's speed must not exceed a specified threshold $v_{\mathrm{save}}$ \textit{until} it stops at the battery charging station, a circular area centered at $r_{\mathrm{w}}$ with radius $d_{\mathrm{w}}$, for $k_{\mathrm{w}}$ consecutive time steps.

In formulating the corresponding STL specification, we initially define two predicate functions along with their respective predicates as follows:
\begin{align*}
    f(x) := v_{\mathrm{save}} - \|v\|, \; \; & \; \; \varphi_1 := (f(x) \geq 0) \\
    g(x) := d_{\mathrm{w}} - \|r - r_w\|, \; \; & \; \; \varphi_2 := (g(x) \geq 0)
\end{align*}
to assess whether the vehicle's speed exceeds the specified threshold and whether the vehicle is positioned at the battery charging station, respectively.

We subsequently define two predicates as follows
\begin{align*}
    \vartheta_1 &:= \bm{G}_{\intv{0}{k_{\mathrm{w}}-1}} \varphi_1 \\
    \vartheta_2 &:= \bm{G}_{\intv{0}{k_{\mathrm{w}}-1}} \varphi_2
\end{align*}
to assess if these specifications are satisfied for $k_{\mathrm{w}}$ consecutive time steps.

Given the scenario's requirement to satisfy the speed specification, \textit{until} the vehicle is positioned at the battery charging station for $k_{\mathrm{w}}$ consecutive time steps, the resultant specification is formulated using D-GMSR as follows:
\begin{equation} \label{eq:qf_stl}
    \Gamma^{ \vartheta_1 \bm{U}_{[0 : K -( k_{\mathrm{w}}-1) ]} \vartheta_2 }_{\bm{\epsilon}, \bm{p}, \bm{w}} (x, 1)
\end{equation}
The construction of the specification using D-GMSR is presented in Appendix \ref{app:qf}.

The resulting problem is presented in (\ref{eq:ct_opt_cont}) where Table \ref{tab:dbl-int-param} shows the system and simulation parameters.
\begin{empheq}[box=\fbox]{equation}
    \begin{aligned} \label{eq:ct_opt_cont}
    \underset{x_k, u_k}{\text{maximize}} \quad & 
        \eqref{eq:qf_stl} \quad
        \text{subject to} \quad 
        & \eqref{eq:quad_dyn_dt}, \eqref{eq:quad_bc}, \eqref{eq:quad_cons}
        \end{aligned}
\end{empheq}

\begin{table}[!htpb]
\centering
\caption{}\label{tab:dbl-int-param}
{\renewcommand{\arraystretch}{1.1}
\begin{tabular}{l|l}
\hline
Parameter & Value\\\hline\\[-0.3cm]
$t_f$, $K$ & $20$ s, $21$ \\
$r_{\mathrm{i}}$, $r_{\mathrm{f}}$ & $(0,0,0)$ m, $(8,0,0)$ m \\
$v_{\mathrm{i}}$, $v_{\mathrm{f}}$ & $0$ m s$^{-1}$, $0$ m s$^{-1}$ \\
$T_{\mathrm{max}}$, $m$ &$10.3$ N, $1.0$ kg \\
$v_{\mathrm{max}}$, $\theta_{\mathrm{max}}$ &  $2$ m s$^{-1}$, $10$ deg \\
$v_{\mathrm{save}}$, $k_{\mathrm{w}}$ &  $1$ m s$^{-1}$, 11 \\
$r_{\mathrm{w}}$, $d_{\mathrm{w}}$ & $(1.25, 2, 2)$ m, $0.2$ m \\
$g$, $c_{\mathrm{d}}$ & $(0, 0, -9.806)$ m s$^{-2}$, $0.5$ m$^{-1}$ \\
$\epsilon_i \in \bm{\epsilon}$, $p_i \in \bm{p}$, $w_i \in \bm{w}$ & $10^{-8}$, $1$, $\bm{1}$ for all $i \in \mathbb{Z}_{++}$\\
\hline
\end{tabular}}
\end{table}

\begin{figure*}[!htpb]
\centering
\begin{subfigure}[b]{1.08\columnwidth}
\includegraphics[width=\linewidth]{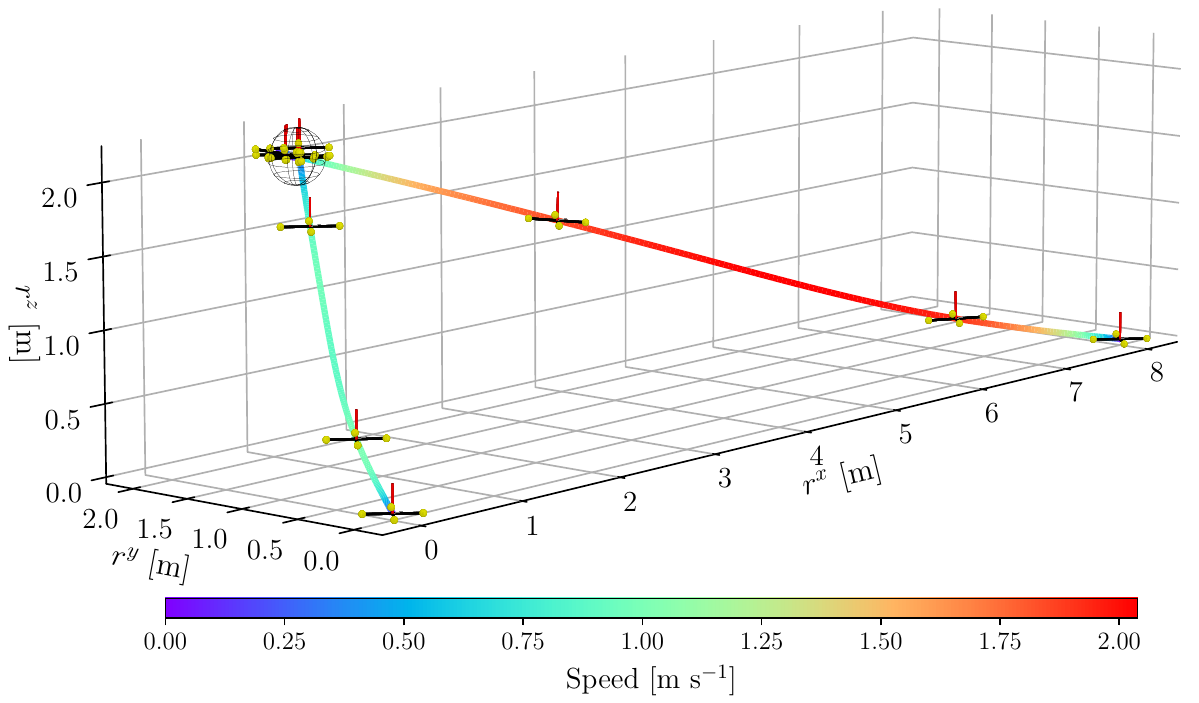}
\label{fig:qf_pos}
\end{subfigure}
\begin{subfigure}[b]{1.01\columnwidth}
\includegraphics[width=\linewidth]{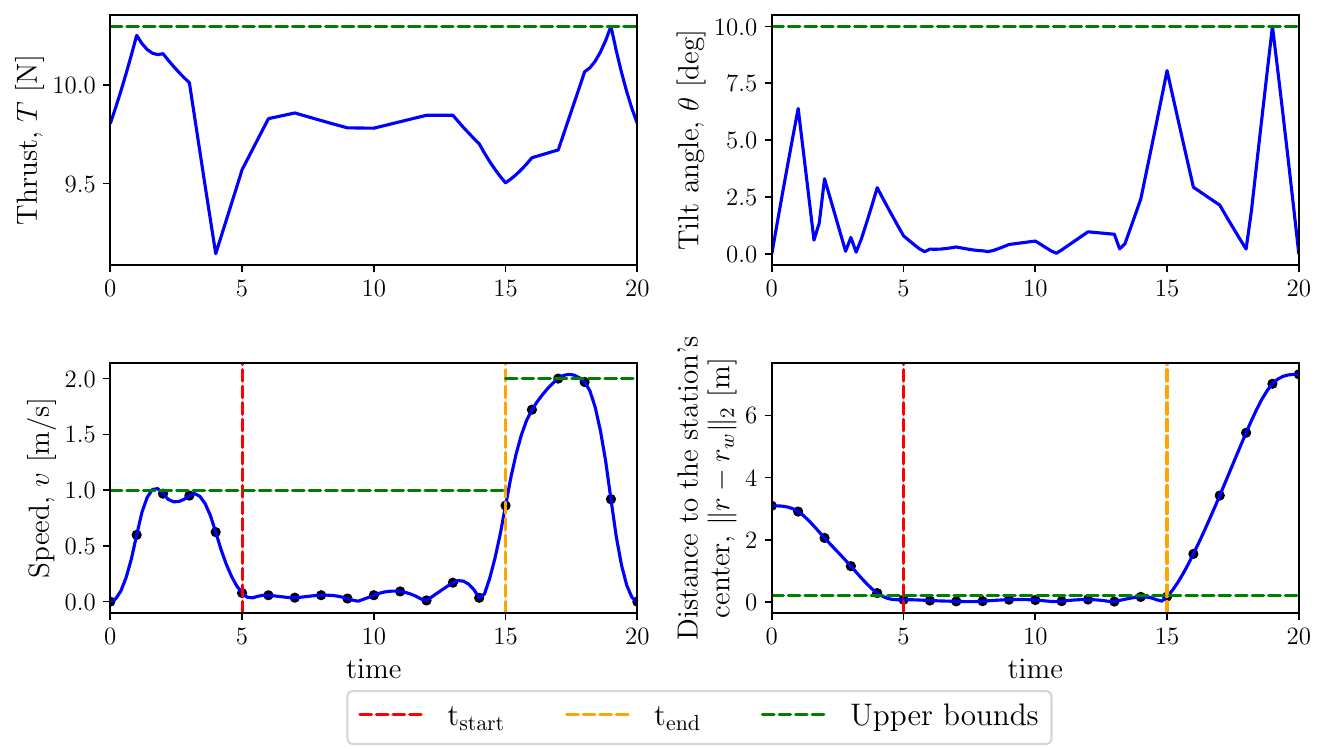}    
\label{fig:qf_oth}
\end{subfigure}
\caption{Quadrotor flight with a battery charging station specification.}
\label{fig:qf_Num_sim}
\end{figure*}

Since \eqref{eq:quad_dyn_dt} is a nonconvex constraint, it is penalized with $\ell_1$ norm, and the resulting optimization problem is solved via the prox-linear method. The result, where the quadrotor executes the control input solution provided by the optimization algorithm, is presented in Figure \ref{fig:qf_Num_sim}.
While the constraints on thrust, tilt angle, and speed are satisfied at the node points, inter-sample constraint violations discussed in Remark \ref{rem:iscv} are observed on the speed constraint.
In conjunction with the reformulation technique proposed in \cite{ctcs2024}, continuous-time formulation of the STL specifications is required to ensure the continuous-time satisfaction of both STL specifications and other constraints.
The quadrotor enters the battery charging station at the $5$th time step and stays there \textit{until} the $15$th time step, ensuring that its speed remains below the threshold throughout this period. To maximize the robustness value, the node points that are in the battery charging station converge through the center of the station. Also, it is observed that the vehicle's speed remains significantly lower than the value of $v_{\textrm{save}}$ \textit{until} it exits the charging station.

\vspace{-0.15cm}
\subsection{Autonomous Rocket landing} \label{sec:rl}
\vspace{-0.15cm}
We address a powered descent guidance problem, specifically, the rocket landing scenario presented in \cite{kamath2023seco}, by incorporating implication specifications, i.e., state-triggered constraints (STCs). 
In this scenario, the rocket must execute a precision landing, involving a maneuver to transition from a horizontal to vertical orientation, and reduce the number of operating engines from $3$ to $1$, similar to Starship \cite{starship2021}.
While the constraint for decreasing the number of operating engines is triggered by speed, multiple constraints related to the landing's safety are triggered by altitude.
\subsubsection{Vehicle dynamics}
The state vector of the vehicle is defined as 
\vspace{-0.15cm}
\begin{equation*}
    x = (m, r, v, \theta, \omega, \delta) \in \mathbb{R}^8
\end{equation*}
where $m \in \mathbb{R}^1_{++}$, $r := (r^x, r^z) \in \mathbb{R}^2$, $v := (v^x, v^z) \in \mathbb{R}^2$, $\theta \in \mathbb{R}^1$, $\omega \in \mathbb{R}^1$ and $\delta \in \mathbb{R}^1$ represent the mass, position, velocity, tilt angle, angular velocity of the vehicle, and gimbal angle of the vehicle's engine, respectively. 
The control input of the vehicle is defined as
\vspace{-0.25cm}
\begin{align*}
    u = (T, \dot{\delta}) \in \mathbb{R}^4
\end{align*}
where $ T := (T^1, T^2, T^3) \in \mathbb{R}_{+}^3 $, with $ T^i \in \mathbb{R}_{+}$ representing the thrust magnitude of the $i$th engine, and $ \dot{\delta} \in \mathbb{R}^1 $denotes the gimbal rate of the vehicle's engine.
The dynamic equation of the vehicle is defined as
    \begin{align} \label{eq:dyn_eq_r}
        \Dot{x}(t) &=  f(x(t), u(t))\\ 
        &= 
        \begin{bmatrix}
            - \alpha \bm{1}^T T(t)  \\
            v \\
            g_{\mathcal{I}} + \frac{1}{m(t)} (F_{\mathcal{I}}(t) + D_{\mathcal{I}}(t)) \\
            \omega \\
            \frac{1}{J_{\mathcal{B}}(t)} ( -l_{\mathrm{cm}} \bm{1}^T T(t) \sin (\delta(t)) - l_{\mathrm{cp}} e_2^{\top} D_{\mathcal{B}}(t) ) \\
            \dot{\delta}
        \end{bmatrix} \nonumber
    \end{align}
    where 
    \begin{align*}
        \alpha &:= \frac{1}{I_{\mathrm{sp}} g_0} \\
        g_{\mathcal{I}} &:= (0, g_0) \\
        F_{\mathcal{I}}(t) &:= \bm{1}^T T(t) \begin{bmatrix}
        \cos ( \theta(t) + \delta(t) ) \\ - \sin ( \theta(t) + \delta(t) )
        \end{bmatrix}\\
        D_{\mathcal{I}}(t) &:= \mathcal{R}_{\mathcal{I} \leftarrow \mathcal{B}}(t) D_{\mathcal{B}}(t) \\
        \mathcal{R}_{\mathcal{I} \leftarrow \mathcal{B}}(t) &:= \begin{bmatrix}
        \cos ( \theta(t) ) & \sin ( \theta(t) ) \\ - \sin ( \theta(t) ) & \cos ( \theta(t) )
        \end{bmatrix} \\
        D_{\mathcal{B}}(t) &:= -\frac{1}{2} \rho_{\mathrm{air}} S_{\mathrm{area}} \| v(t) \|_2 C_{\mathrm{aero}} \mathcal{R}_{\mathcal{I} \leftarrow \mathcal{B}}^{\top}(t) v(t) \\
        J_{\mathcal{B}}(t) &:= m(t) \bigg(\frac{l_{\mathrm{r}}^2}{4} + \frac{l_{\mathrm{h}}^2}{12}\bigg) 
    \end{align*}
    Here, 
    $\alpha \in \mathbb{R}_{++}$ is  
    the thrust-specific fuel consumption, 
    $I_{\mathrm{sp}} \in \mathbb{R}_{++}$ is the specific impulse of the rocket engine, 
    $g_0 \in \mathbb{R}_{++}$ is the average gravitational acceleration of the planet, 
    $F_{\mathcal{I}}(t)$ and $D_{\mathcal{I}}(t)$ are the thrust and drag forces on the vehicle with respect to the inertial frame, respectively,
    $\mathcal{R}_{\mathcal{I} \leftarrow \mathcal{B}}(t) \in \mathrm{SO}(2)$ is the rotation matrix that maps coordinates in the body frame to the inertial frame, 
    $D_{\mathcal{B}}(t)$ is the drag force on the vehicle with respect to the body frame,
    $\rho_{\mathrm{air}} \in \mathbb{R}_{++}$ is the atmospheric density, 
    $S_{\mathrm{area}} \in \mathbb{R}_{++}$ is the reference area, 
    $C_{\mathrm{aero}} := \mathrm{diag} \{c_x, c_z\} $ is the aerodynamic coefficient matrix, where $c_x$, $c_z \in \mathbb{R}_{++}$ are the aerodynamic coefficients along the body $x$ and $z$ axes, respectively,
    $J_{\mathcal{B}}(t) \in \mathbb{R}_{++}$ is the moment of inertia of the vehicle about the body $y$ axis,
    $l_{\mathrm{r}} \in \mathbb{R}_{++}$ and $l_{\mathrm{h}} \in \mathbb{R}_{++}$ are the radius and height of the fuselage, respectively,
    $l_{\mathrm{cm}} \in \mathbb{R}_{++}$ is the thrust moment-arm (distance between the vehicle mass-center and the engine gimbal hinge point),
    and $l_{\mathrm{cp}} \in \mathbb{R}_{+}$ is the aerodynamic moment-arm (the distance between the vehicle mass-center and the center of pressure).

    We discretize the continuous-time dynamics using the multiple-shooting approach. The time is discretized over a finite grid as $0 = t_1 < t_2 < \dots < t_K = t_f$, $x_k$ and $u_k$ represent the state and control input at the node point $t_k$ for $k \in \intv{1}{K}$. The control input $u$ is parameterized via first-order hold and the continuous-time dynamics \eqref{eq:quad_dyn_eq} is discretized using multiple-shooting as follows:
    \begin{align} \label{eq:rl_dyn_dt}
        x_{k+1} = F_k(x_k, u_k, u_{k+1})
    \end{align}
    for $k \in \intv{1}{K-1}$. See \cite{ctcs2024, nmpc2024} for further details.

    The boundary conditions are specified as follows:
    \begin{subequations} \label{eq:rl_bc}
    \begin{gather}
        m_1 = m_{\mathrm{i}}, \; m_K \geq m_{\mathrm{dry}}, \; r_1 = r_{\mathrm{i}}, \; r_K = r_{\mathrm{f}} \tag{\theequation a-d} \\
        v_1 = v_{\mathrm{i}}, \; v_K = v_{\mathrm{f}}, \; \theta_1 = \theta_{\mathrm{i}}, \; \theta_K = \theta_{\mathrm{f}} \tag{\theequation e-h} \\
        \!\!\!\!\!\!\!\!\!\!\!\!\!\!\!\!\!\!\!\!\! \omega_1 = \omega_{\mathrm{i}}, \; \omega_K = \omega_{\mathrm{f}}, \; \delta_1 = \delta_{\mathrm{i}} \tag{\theequation i-k} 
    \end{gather}
    \end{subequations}
    where $m_{\mathrm{dry}}$ is the dry mass of the vehicle.
    
    The altitude, gimbal angle, thrust, and gimbal rate are constrained as follows:
    \begin{subequations} \label{eq:rl_cons}
    \begin{align}
        r_{\mathrm{min}}^x &\leq r^x_k \\
        \| \delta_k \| &\leq \delta_{\mathrm{max}} \\
        T_{\mathrm{min}} &\leq T_k^1 \leq T_{\mathrm{max}} \\
        0 &\leq T_k^2 \leq T_{\mathrm{max}} \\
        0 &\leq T_k^3 \leq T_{\mathrm{max}} \\
        \| \dot{\delta_k} \| &\leq \dot{\delta}_{\mathrm{max}}
    \end{align}
    \end{subequations}
    for $k \in \intv{1}{K-1}$, where $r_{\mathrm{min}}^x$ is the minimum altitude.

    \subsubsection{STL specifications} We impose $14$ different implication specifications i.e. STC, where constraints are triggered by the speed and altitude, as follows:
    \begin{align*}
        \varphi_{i} :=& \;  ( \| v_{k'} \|_2 > v_{\mathrm{trig}} ) \! \implies \! ( T_{\mathrm{min}} \leq T_k^{i+1}) \text{  for $i \!=\! 1, 2$} \\
        \varphi_{i+2} :=& \;  ( \| v_{k'} \|_2 < v_{\mathrm{trig}} ) \! \implies \! ( T_k^{i+1} = 0) \text{  for $i \!=\! 1, 2$} \\
        \varphi_5 :=& \;  ( \| v_{k'} \|_2 < v_{\mathrm{trig}} ) \! \implies \! ( \| v_k \|_2 + \epsilon_{\mathrm{v}} \leq v_{\mathrm{trig}} )\\
        \varphi_6 :=& \;  ( r_{k'}^x <r_{\mathrm{trig}}^x ) \! \implies \! ( \| v_k \| \leq v_{\mathrm{stc}} )\\
        \varphi_7 :=& \;  ( r_{k'}^x <r_{\mathrm{trig}}^x ) \! \implies \! ( \theta_k \leq \theta_{\mathrm{stc}} )\\
        \varphi_8 :=& \;  ( r_{k'}^x <r_{\mathrm{trig}}^x ) \! \implies \! ( \theta_k \geq -\theta_{\mathrm{stc}} )\\
        \varphi_9 :=& \;  ( r_{k'}^x <r_{\mathrm{trig}}^x ) \! \implies \! ( \omega_k \leq \omega_{\mathrm{stc}} )\\
        \varphi_{10} :=& \;  ( r_{k'}^x <r_{\mathrm{trig}}^x ) \! \implies \! ( \omega_k \geq -\omega_{\mathrm{stc}} )\\
        \varphi_{11} :=& \;  ( r_{k'}^x <r_{\mathrm{trig}}^x ) \! \implies \! ( \delta_k \leq \delta_{\mathrm{stc}} )\\
        \varphi_{12} :=& \;  ( r_{k'}^x <r_{\mathrm{trig}}^x ) \! \implies \! ( \delta_k \geq -\delta_{\mathrm{stc}} )\\
        \varphi_{13} :=& \; ( r_{k'}^x <r_{\mathrm{trig}}^x ) \! \implies \! ( r_{z_k} \leq r^x_k \tan(\gamma_{\mathrm{stc}})
        )\\
        \varphi_{14} :=& \;  ( r_{k'}^x <r_{\mathrm{trig}}^x ) \! \implies \! ( r_{z_k} \geq - r^x_k \tan(\gamma_{\mathrm{stc}}) )
    \end{align*}
    $\forall k \in \intv{k'}{K}$, where $v_{\mathrm{trig}}$ and $r_{\mathrm{trig}}^x$ are the velocity and altitude that trigger the constraints, respectively, $\gamma_{\mathrm{max}}$ is the maximum glideslope angle, and $\epsilon_{\mathrm{v}}$ is an arbitrarily small positive number. 
    
    All specifications are formulated via D-GMSR and the following constraint is imposed on the optimization problem:
    \begin{equation} \label{eq:rl_stl}
        \Gamma_{\bm{\epsilon}, \bm{p}, \bm{w}}^{\varphi_i} (x,1) \geq 0 \text{ for } i = 1,2,\dots, 14 
    \end{equation}
    The objective function of this problem is to minimize the fuel consumption of the rocket.
    The resulting optimization problem is presented in (\ref{eq:rl_ct_opt_cont}) where Table \ref{tab:rl-int-param} shows the system and simulation parameters.
    \begin{empheq}[box=\fbox]{equation}
        \begin{aligned} \label{eq:rl_ct_opt_cont}
        \underset{x_k, u_k}{\text{maximize}} \;\; 
            m_K \;\;
            \text{subject to} \;
            \eqref{eq:rl_dyn_dt}, \eqref{eq:rl_bc}, \eqref{eq:rl_cons},
            \eqref{eq:rl_stl}
        \end{aligned}
    \end{empheq}
    \begin{figure*}[!htpb]
    \centering
    \begin{subfigure}[b]{0.682\columnwidth}
    \includegraphics[width=\linewidth]{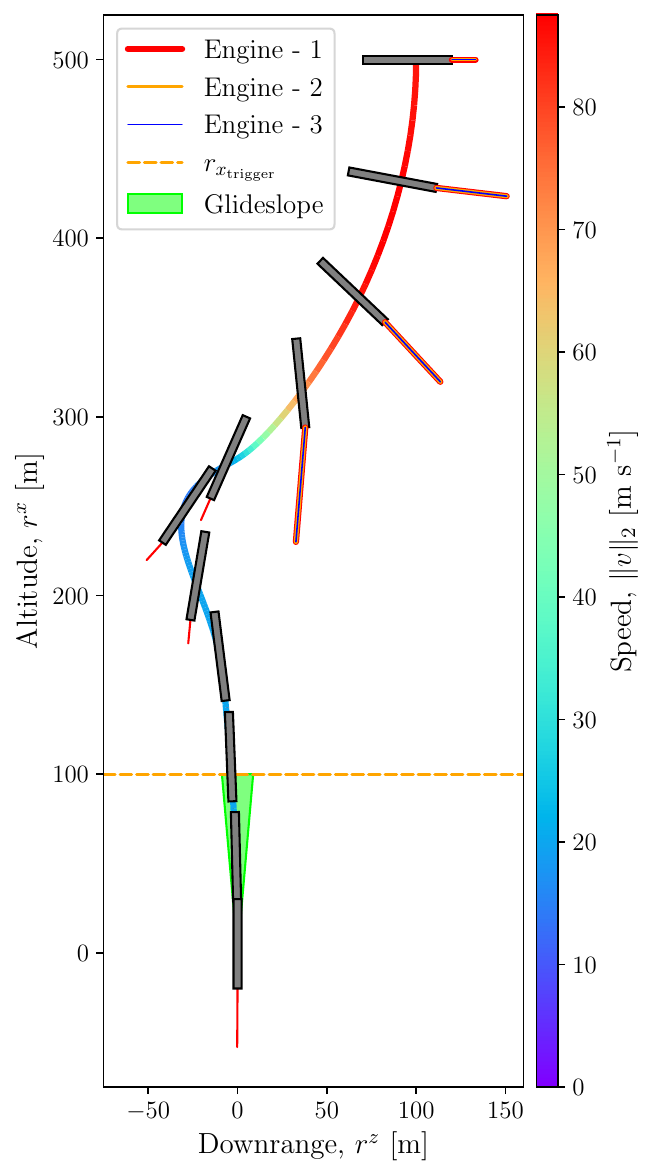}
    \label{fig:rl_pos}
    \end{subfigure}
    \begin{subfigure}[b]{1.409\columnwidth}
    \includegraphics[width=\linewidth]{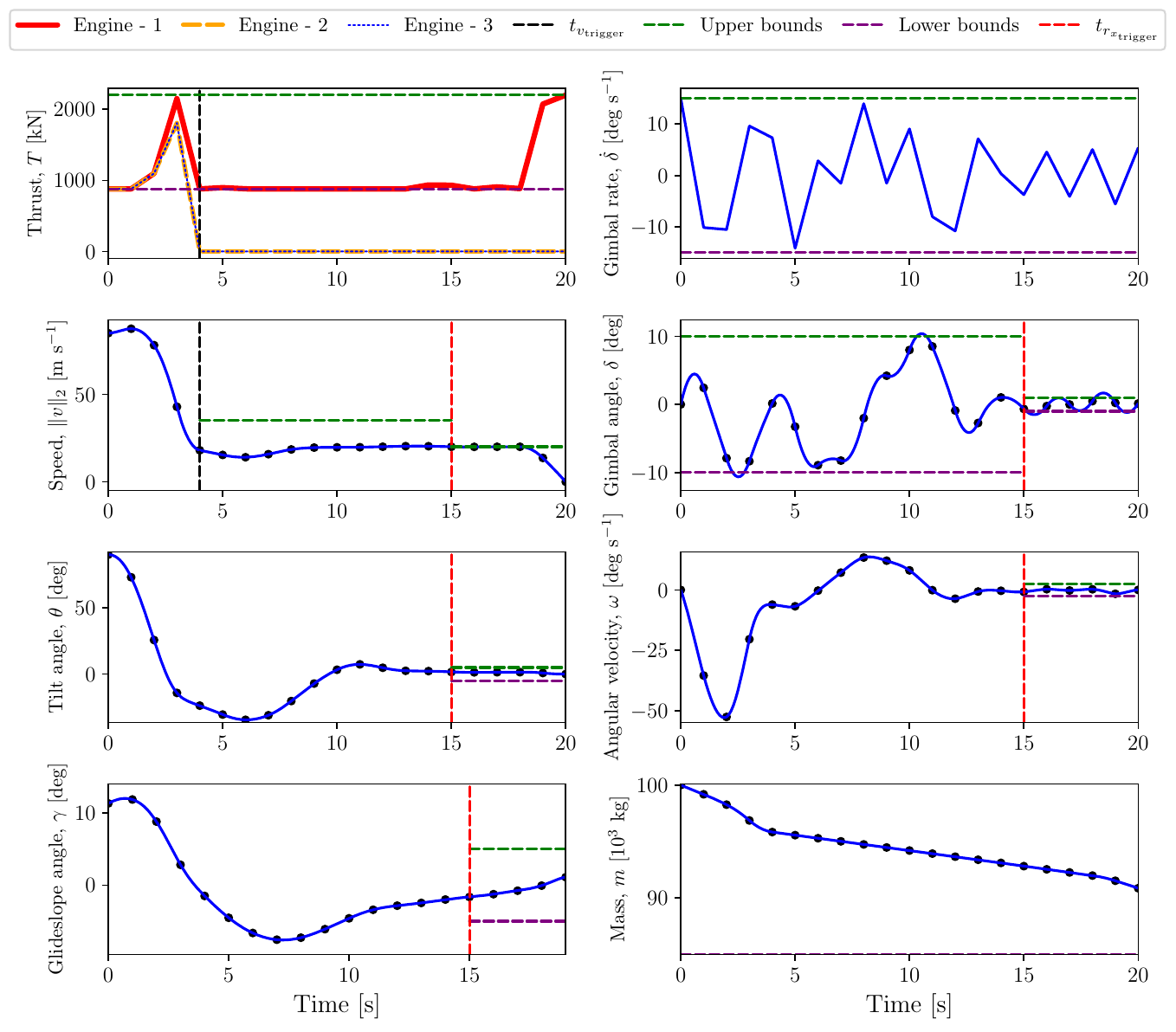}    
    \label{fig:rl_oth}
    \end{subfigure}
    \caption{Rocket landing with state-triggered constraints.}
    \label{fig:rl_num_sim}
    \end{figure*}
    \begin{table}[!htpb]
    \centering
    \caption{}\label{tab:rl-int-param}
    {\renewcommand{\arraystretch}{1.1}
    \begin{tabular}{l|l}
    \hline
    Parameter & Value\\\hline\\[-0.3cm]
    $t_f$, $K$ & $20$ s, $21$ \\
    $I_\mathrm{sp}$, $g_0$ & $330$ s, $9.806$ m s$^{-1}$ \\
    $l_{\mathrm{r}}$, $l_{\mathrm{h}}$ & $4.5$ m, $50$ m \\
    $l_{\mathrm{cm}}$, $l_{\mathrm{cp}}$ & $20$ m, $10$ m \\
    $\rho_{\mathrm{air}}$, $S_{\mathrm{area}}$ & $1.225$ kg m$^{-3}$, $545$ m$^2$ \\
    $c_{x}$, $c_{z}$ & $0.0522$, $0.4068$ \\
    $m_{\mathrm{i}}$, $m_{\mathrm{dry}}$ & $100000$ kg, $85000$ kg \\
    $r_{\mathrm{i}}$, $r_{\mathrm{f}}$ & $(500,100)$ m, $(0,0)$ m \\
    $v_{\mathrm{i}}$, $v_{\mathrm{f}}$ & $(-85, 0)$ m s$^{-1}$, $(0,0)$ m s$^{-1}$ \\
    $\theta_{\mathrm{i}}$, $\theta_{\mathrm{f}}$ & $90$ deg, $0$ deg \\
    $\omega_{\mathrm{i}}$, $\omega_{\mathrm{f}}$,
    $\delta_{\mathrm{i}}$ & $0$ deg s$^{-1}$, $0$ deg s$^{-1}$, $0$ deg \\
    $r_{\mathrm{min}}^x$, $\delta_{\mathrm{max}}$ & $0$ m $10$ deg\\
    $T_{\mathrm{max}}$, $T_{\mathrm{min}}$, $\dot{\delta}_{\mathrm{max}}$ &  $2200$ kN, $880$ kN, $15$ deg s$^{-1}$  \\
    $v_{\mathrm{trig}}$, $r_{\mathrm{trig}}^x$, $v_{\mathrm{stc}}$ & $35$ m s$^{-1}$, $100$ m, $20$ m s$^{-1}$ \\
    $\theta_{\mathrm{stc}}$, $\omega_{\mathrm{stc}}$, $\delta_{\mathrm{stc}}$, 
    $\gamma_{\mathrm{stc}}$ 
    & $5$ deg, $2.5$ deg s$^{-1}$, $1$ deg, $5$ deg\\
    $\epsilon_i \in \bm{\epsilon}$, $p_i \in \bm{p}$, $w_i \in \bm{w}$ & $10^{-8}$, $1$, $\bm{1}$ for all $i \in \mathbb{Z}_{++}$\\
    \hline
    \end{tabular}}
    \end{table}

    \vspace{-0.3cm}
    Since \eqref{eq:rl_dyn_dt} and \eqref{eq:rl_stl} are nonconvex constraints, they are penalized with $\ell_1$ norm, and the resulting optimization problem is solved via the prox-linear method. 
    The result, where the rocket executes the control input solution provided by the optimization algorithm, is presented in Figure \ref{fig:rl_num_sim}.
    While the constraints are satisfied at the node points, inter-sample constraint violations discussed in Remark \ref{rem:iscv} are observed on the gimbal angle constraint. 
    In conjunction with the reformulation technique proposed in \cite{ctcs2024}, continuous-time formulation of the STL specifications is required to ensure the continuous-time satisfaction of both STL specifications and other constrains.
    During precision landing, the vehicle performs a maneuver to transition from horizontal to vertical orientation, akin to Starship \cite{starship2021}, and it satisfies all constraints triggered by speed and altitude.

\vspace{-0.15cm}
\section{Conclusion and Future work} \label{sec:conc}
\vspace{-0.15cm}

This paper introduces a $\mathcal{C}^1$-smooth, sound, and complete robustness measure for STL specifications. The robustness measure demonstrates desirable gradient properties and therefore addresses the locality and masking problems, which are critical for numerical optimization.
The performance of the robustness measure is demonstrated on real-world trajectory optimization problems for quadrotor flight and autonomous rocket landing solved within an SCP framework.

The continuous-time modeling of STL specifications and a successive convexification based solution method to ensure the continuous-time satisfaction of STL specifications will be the subject of future studies. 
Additionally, potential future work may involve the design of optimization algorithms that exploit the inherent structure of STL specifications to achieve faster convergence.
\begin{ack}
The authors gratefully acknowledge Fabio Spada and Abhinav G. Kamath for helpful discussions and their feedback on the initial drafts of the paper.
\end{ack}  
{\RaggedRight%
\bibliographystyle{unsrturl}      
\bibliography{references}%

\begin{thebibliography}{10}

\bibitem{baier2008principles}
Christel Baier and Joost-Pieter Katoen.
\newblock {\em Principles of Model Checking}.
\newblock MIT press, 2008.
\newblock URL:
  \url{https://mitpress.mit.edu/9780262026499/principles-of-model-checking/}.

\bibitem{sahin2020autonomous}
Yunus~Emre Sahin, Rien Quirynen, and Stefano Di~Cairano.
\newblock Autonomous vehicle decision-making and monitoring based on signal
  temporal logic and mixed-integer programming.
\newblock In {\em 2020 American Control Conference (ACC)}, pages 454--459.
  IEEE, 2020.
\newblock URL: \url{https://doi.org/10.23919/ACC45564.2020.9147917}.

\bibitem{kurtz2020trajectory}
Vince Kurtz and Hai Lin.
\newblock Trajectory optimization for high-dimensional nonlinear systems under
  {STL} specifications.
\newblock {\em IEEE Control Systems Letters}, 5(4):1429--1434, 2020.
\newblock URL: \url{https://doi.org/10.1109/LCSYS.2020.3038640}.

\bibitem{djeumou2022probabilistic}
Franck Djeumou, Zhe Xu, Murat Cubuktepe, and Ufuk Topcu.
\newblock Probabilistic control of heterogeneous swarms subject to graph
  temporal logic specifications: A decentralized and scalable approach.
\newblock {\em IEEE Transactions on Automatic Control}, 68(4):2245--2260, 2022.
\newblock URL: \url{https://doi.org/10.1109/TAC.2022.3176797}.

\bibitem{sahin2019multirobot}
Yunus~Emre Sahin, Petter Nilsson, and Necmiye Ozay.
\newblock Multirobot coordination with counting temporal logics.
\newblock {\em IEEE Transactions on Robotics}, 36(4):1189--1206, 2019.
\newblock URL: \url{https://doi.org/10.1109/TRO.2019.2957669}.

\bibitem{buyukkocak2021planning}
Ali~Tevfik Buyukkocak, Derya Aksaray, and Yasin Yaz{\i}c{\i}o{\u{g}}lu.
\newblock Planning of heterogeneous multi-agent systems under signal temporal
  logic specifications with integral predicates.
\newblock {\em IEEE Robotics and Automation Letters}, 6(2):1375--1382, 2021.
\newblock URL: \url{https://doi.org/10.1109/LRA.2021.3057049}.

\bibitem{pant2018fly}
Yash~Vardhan Pant, Houssam Abbas, Rhudii~A Quaye, and Rahul Mangharam.
\newblock Fly-by-logic: Control of multi-drone fleets with temporal logic
  objectives.
\newblock In {\em 2018 ACM/IEEE 9th International Conference on Cyber-Physical
  Systems (ICCPS)}, pages 186--197. IEEE, 2018.
\newblock URL: \url{https://doi.org/10.1109/ICCPS.2018.00026}.

\bibitem{bacspinar2019mission}
Bari{\c{s}} Ba{\c{s}}pinar, Hamsa Balakrishnan, and Emre Koyuncu.
\newblock Mission planning and control of multi-aircraft systems with signal
  temporal logic specifications.
\newblock {\em IEEE Access}, 7:155941--155950, 2019.
\newblock URL: \url{https://doi.org/10.1109/ACCESS.2019.2949426}.

\bibitem{ma2020stlnet}
Meiyi Ma, Ji~Gao, Lu~Feng, and John Stankovic.
\newblock {STLnet}: Signal temporal logic enforced multivariate recurrent
  neural networks.
\newblock {\em Advances in Neural Information Processing Systems},
  33:14604--14614, 2020.
\newblock URL: \url{https://dl.acm.org/doi/abs/10.5555/3495724.3496948}.

\bibitem{aksaray2016q}
Derya Aksaray, Austin Jones, Zhaodan Kong, Mac Schwager, and Calin Belta.
\newblock Q-learning for robust satisfaction of signal temporal logic
  specifications.
\newblock In {\em 2016 IEEE 55th Conference on Decision and Control (CDC)},
  pages 6565--6570. IEEE, 2016.
\newblock URL: \url{https://doi.org/10.1109/CDC.2016.7799279}.

\bibitem{alshiekh2018safe}
Mohammed Alshiekh, Roderick Bloem, R{\"u}diger Ehlers, Bettina K{\"o}nighofer,
  Scott Niekum, and Ufuk Topcu.
\newblock Safe reinforcement learning via shielding.
\newblock In {\em Proceedings of the AAAI Conference on Artificial
  Intelligence}, volume~32, 2018.
\newblock URL: \url{https://doi.org/10.1609/AAAI.V32I1.11797}.

\bibitem{leung2023backpropagation}
Karen Leung, Nikos Ar{\'e}chiga, and Marco Pavone.
\newblock Backpropagation through signal temporal logic specifications:
  Infusing logical structure into gradient-based methods.
\newblock {\em The International Journal of Robotics Research}, 42(6):356--370,
  2023.
\newblock URL: \url{https://doi.org/10.1177/02783649221082115}.

\bibitem{pnueli1977temporal}
Amir Pnueli.
\newblock The temporal logic of programs.
\newblock In {\em 18th Annual Symposium on Foundations of Computer Science
  (SFCS 1977)}, pages 46--57. IEEE, 1977.
\newblock URL: \url{https://doi.org/10.1109/SFCS.1977.32}.

\bibitem{alur1994theory}
Rajeev Alur and David~L Dill.
\newblock A theory of timed automata.
\newblock {\em Theoretical Computer Science}, 126(2):183--235, 1994.
\newblock URL: \url{https://doi.org/10.1016/0304-3975(94)90010-8}.

\bibitem{wongpiromsarn2012receding}
Tichakorn Wongpiromsarn, Ufuk Topcu, and Richard~M Murray.
\newblock Receding horizon temporal logic planning.
\newblock {\em IEEE Transactions on Automatic Control}, 57(11):2817--2830,
  2012.
\newblock URL: \url{https://doi.org/10.1109/TAC.2012.2195811}.

\bibitem{wongpiromsarn2011tulip}
Tichakorn Wongpiromsarn, Ufuk Topcu, Necmiye Ozay, Huan Xu, and Richard~M
  Murray.
\newblock {TuLiP}: a software toolbox for receding horizon temporal logic
  planning.
\newblock In {\em Proceedings of the 14th International Conference on Hybrid
  Systems: Computation and Control}, pages 313--314, 2011.
\newblock URL: \url{https://doi.org/10.1145/1967701.1967747}.

\bibitem{karaman2008optimal}
Sertac Karaman, Ricardo~G Sanfelice, and Emilio Frazzoli.
\newblock Optimal control of mixed logical dynamical systems with linear
  temporal logic specifications.
\newblock In {\em 2008 47th IEEE Conference on Decision and Control}, pages
  2117--2122. IEEE, 2008.
\newblock URL: \url{https://doi.org/10.1109/CDC.2008.4739370}.

\bibitem{koymans1990specifying}
Ron Koymans.
\newblock Specifying real-time properties with metric temporal logic.
\newblock {\em Real-Time Systems}, 2(4):255--299, 1990.
\newblock URL: \url{https://doi.org/10.1007/BF01995674}.

\bibitem{alur1996benefits}
Rajeev Alur, Tom{\'a}s Feder, and Thomas~A Henzinger.
\newblock The benefits of relaxing punctuality.
\newblock {\em Journal of the ACM (JACM)}, 43(1):116--146, 1996.
\newblock URL: \url{https://doi.org/10.1145/227595.227602}.

\bibitem{maler2004monitoring}
Oded Maler and Dejan Nickovic.
\newblock Monitoring temporal properties of continuous signals.
\newblock In {\em International Symposium on Formal Techniques in Real-Time and
  Fault-Tolerant Systems}, pages 152--166. Springer, 2004.
\newblock URL: \url{https://doi.org/10.1007/978-3-540-30206-3_12}.

\bibitem{fainekos2009robustness}
Georgios~E Fainekos and George~J Pappas.
\newblock Robustness of temporal logic specifications for continuous-time
  signals.
\newblock {\em Theoretical Computer Science}, 410(42):4262--4291, 2009.
\newblock URL: \url{https://doi.org/10.1016/j.tcs.2009.06.021}.

\bibitem{donze2010robust}
Alexandre Donz{\'e} and Oded Maler.
\newblock Robust satisfaction of temporal logic over real-valued signals.
\newblock In {\em International Conference on Formal Modeling and Analysis of
  Timed Systems}, pages 92--106. Springer, 2010.
\newblock URL: \url{https://doi.org/10.1007/978-3-642-15297-9_9}.

\bibitem{lindemann2018control}
Lars Lindemann and Dimos~V Dimarogonas.
\newblock Control barrier functions for signal temporal logic tasks.
\newblock {\em IEEE Control Systems Letters}, 3(1):96--101, 2018.
\newblock URL: \url{https://doi.org/10.1109/LCSYS.2018.2853182}.

\bibitem{belta2019formal}
Calin Belta and Sadra Sadraddini.
\newblock Formal methods for control synthesis: An optimization perspective.
\newblock {\em Annual Review of Control, Robotics, and Autonomous Systems},
  2:115--140, 2019.
\newblock URL: \url{https://doi.org/10.1146/annurev-control-053018-023717}.

\bibitem{raman2014model}
Vasumathi Raman, Alexandre Donz{\'e}, Mehdi Maasoumy, Richard~M Murray, Alberto
  Sangiovanni-Vincentelli, and Sanjit~A Seshia.
\newblock Model predictive control with signal temporal logic specifications.
\newblock In {\em 53rd IEEE Conference on Decision and Control}, pages 81--87.
  IEEE, 2014.
\newblock URL: \url{https://doi.org/10.1109/CDC.2014.7039363}.

\bibitem{sadraddini2015robust}
Sadra Sadraddini and Calin Belta.
\newblock Robust temporal logic model predictive control.
\newblock In {\em 2015 53rd Annual Allerton Conference on Communication,
  Control, and Computing (Allerton)}, pages 772--779. IEEE, 2015.
\newblock URL: \url{https://doi.org/10.1109/ALLERTON.2015.7447084}.

\bibitem{rodionova2021time}
Al{\"e}na Rodionova, Lars Lindemann, Manfred Morari, and George~J Pappas.
\newblock Time-robust control for {STL} specifications.
\newblock In {\em 2021 60th IEEE Conference on Decision and Control (CDC)},
  pages 572--579. IEEE, 2021.
\newblock URL: \url{https://doi.org/10.1109/CDC45484.2021.9683477}.

\bibitem{rodionova2022combined}
Al{\"e}na Rodionova, Lars Lindemann, Manfred Morari, and George~J Pappas.
\newblock Combined left and right temporal robustness for control under {STL}
  specifications.
\newblock {\em IEEE Control Systems Letters}, 7:619--624, 2022.
\newblock URL: \url{https://doi.org/10.1109/LCSYS.2022.3209928}.

\bibitem{pant2017smooth}
Yash~Vardhan Pant, Houssam Abbas, and Rahul Mangharam.
\newblock Smooth operator: Control using the smooth robustness of temporal
  logic.
\newblock In {\em 2017 IEEE Conference on Control Technology and Applications
  (CCTA)}, pages 1235--1240. IEEE, 2017.
\newblock URL: \url{https://doi.org/10.1109/CCTA.2017.8062628}.

\bibitem{haghighi2019control}
Iman Haghighi, Noushin Mehdipour, Ezio Bartocci, and Calin Belta.
\newblock Control from signal temporal logic specifications with smooth
  cumulative quantitative semantics.
\newblock In {\em 2019 IEEE 58th Conference on Decision and Control (CDC)},
  pages 4361--4366. IEEE, 2019.
\newblock URL: \url{https://doi.org/10.1109/CDC40024.2019.9029429}.

\bibitem{gilpin2020smooth}
Yann Gilpin, Vince Kurtz, and Hai Lin.
\newblock A smooth robustness measure of signal temporal logic for symbolic
  control.
\newblock {\em IEEE Control Systems Letters}, 5(1):241--246, 2020.
\newblock URL: \url{https://doi.org/10.1109/LCSYS.2020.3001875}.

\bibitem{mao2022successive}
Yuanqi Mao, Behcet Acikmese, Pierre-Loic Garoche, and Alexandre Chapoutot.
\newblock Successive convexification for optimal control with signal temporal
  logic specifications.
\newblock In {\em Proceedings of the 25th ACM International Conference on
  Hybrid Systems: Computation and Control}, pages 1--7, 2022.
\newblock URL: \url{https://doi.org/10.1145/3501710.3519518}.

\bibitem{lindemann2019robust}
Lars Lindemann and Dimos~V Dimarogonas.
\newblock Robust control for signal temporal logic specifications using
  discrete average space robustness.
\newblock {\em Automatica}, 101:377--387, 2019.
\newblock URL: \url{https://doi.org/10.1016/j.automatica.2018.12.022}.

\bibitem{mehdipour2019arithmetic}
Noushin Mehdipour, Cristian-Ioan Vasile, and Calin Belta.
\newblock Arithmetic-geometric mean robustness for control from signal temporal
  logic specifications.
\newblock In {\em 2019 American Control Conference (ACC)}, pages 1690--1695.
  IEEE, 2019.
\newblock URL: \url{https://doi.org/10.23919/ACC.2019.8814487}.

\bibitem{mehdipour2020specifying}
Noushin Mehdipour, Cristian-Ioan Vasile, and Calin Belta.
\newblock Specifying user preferences using weighted signal temporal logic.
\newblock {\em IEEE Control Systems Letters}, 5(6):2006--2011, 2020.
\newblock URL: \url{https://doi.org/10.1109/LCSYS.2020.3047362}.

\bibitem{mehdipour2019average}
Noushin Mehdipour, Cristian-Ioan Vasile, and Calin Belta.
\newblock Average-based robustness for continuous-time signal temporal logic.
\newblock In {\em 2019 IEEE 58th Conference on Decision and Control (CDC)},
  pages 5312--5317. IEEE, 2019.
\newblock URL: \url{https://doi.org/10.1109/CDC40024.2019.9029989}.

\bibitem{akazaki2015time}
Takumi Akazaki and Ichiro Hasuo.
\newblock Time robustness in {MTL} and expressivity in hybrid system
  falsification.
\newblock In {\em International Conference on Computer Aided Verification},
  pages 356--374. Springer, 2015.
\newblock URL: \url{https://doi.org/10.1007/978-3-319-21668-3_21}.

\bibitem{mao2016successive}
Yuanqi Mao, Michael Szmuk, and Beh{\c{c}}et A{\c{c}}{\i}kme{\c{s}}e.
\newblock Successive convexification of non-convex optimal control problems and
  its convergence properties.
\newblock In {\em 2016 IEEE 55th Conference on Decision and Control (CDC)},
  pages 3636--3641. IEEE, 2016.
\newblock URL: \url{https://doi.org/10.1109/CDC.2016.7798816}.

\bibitem{malyuta2021advances}
Danylo Malyuta, Yue Yu, Purnanand Elango, and Beh{\c{c}}et
  A{\c{c}}{\i}kme{\c{s}}e.
\newblock Advances in trajectory optimization for space vehicle control.
\newblock {\em Annual Reviews in Control}, 52:282--315, 2021.
\newblock URL: \url{https://doi.org/10.1016/j.arcontrol.2021.04.013}.

\bibitem{malyuta2022convex}
Danylo Malyuta, Taylor~P Reynolds, Michael Szmuk, Thomas Lew, Riccardo Bonalli,
  Marco Pavone, and Beh{\c{c}}et A{\c{c}}{\i}kme{\c{s}}e.
\newblock Convex optimization for trajectory generation: A tutorial on
  generating dynamically feasible trajectories reliably and efficiently.
\newblock {\em IEEE Control Systems Magazine}, 42(5):40--113, 2022.
\newblock URL: \url{https://doi.org/10.1109/MCS.2022.3187542}.

\bibitem{ctcs2024}
Purnanand Elango, Dayou Luo, Abhinav~G. Kamath, Samet Uzun, Taewan Kim, and
  Beh{\c{c}}et A{\c{c}}{\i}kme{\c{s}}e.
\newblock Successive convexification for trajectory optimization with
  continuous-time constraint satisfaction.
\newblock {\em arXiv preprint arXiv:2404.16826}, 2024.
\newblock URL: \url{https://doi.org/10.48550/arXiv.2404.16826}.

\bibitem{malyuta2023fast}
Danylo Malyuta and Beh{\c{c}}et A{\c{c}}{\i}kme{\c{s}}e.
\newblock Fast homotopy for spacecraft rendezvous trajectory optimization with
  discrete logic.
\newblock {\em Journal of Guidance, Control, and Dynamics}, 46(7):1262--1279,
  2023.
\newblock URL: \url{https://doi.org/10.2514/1.G006295}.

\bibitem{shniad1948convexity}
Harold Shniad.
\newblock On the convexity of mean value functions.
\newblock {\em Bulletin of the American Mathematical Society}, 54:770--776,
  1948.
\newblock URL: \url{https://doi.org/10.1090/S0002-9904-1948-09077-2}.

\bibitem{tang2018fuel}
Gao Tang, Fanghua Jiang, and Junfeng Li.
\newblock Fuel-optimal low-thrust trajectory optimization using indirect method
  and successive convex programming.
\newblock {\em IEEE Transactions on Aerospace and Electronic Systems},
  54(4):2053--2066, 2018.
\newblock URL: \url{https://doi.org/10.1109/TAES.2018.2803558}.

\bibitem{spada2023direct}
Fabio Spada, Marco Sagliano, and Francesco Topputo.
\newblock Direct--indirect hybrid strategy for optimal powered descent and
  landing.
\newblock {\em Journal of Spacecraft and Rockets}, 60(6):1787--1804, 2023.
\newblock URL: \url{https://doi.org/10.2514/1.A35650}.

\bibitem{Elango2022}
Purnanand Elango, Abhinav Kamath, Yue Yu, John~M. Carson~III, and Behcet
  Acikmese.
\newblock A customized first-order solver for real-time powered-descent
  guidance.
\newblock In {\em AIAA SciTech 2022 Forum}. American Institute of Aeronautics
  and Astronautics, January 2022.
\newblock URL: \url{https://doi.org/10.2514/6.2022-0951}.

\bibitem{Kamath2023}
Abhinav~G. Kamath, Purnanand Elango, Skye Mceowen, Yue Yu, John~M. Carson~III,
  Mehran Mesbahi, and Behcet Acikmese.
\newblock Customized real-time first-order methods for onboard dual
  quaternion-based 6-{DoF} powered-descent guidance.
\newblock In {\em AIAA SciTech 2023 Forum}. American Institute of Aeronautics
  and Astronautics, January 2023.
\newblock URL: \url{https://doi.org/10.2514/6.2023-2003}.

\bibitem{szmuk2019real}
Michael Szmuk, Danylo Malyuta, Taylor~P Reynolds, Margaret~Skye Mceowen, and
  Beh{\c{c}}et A{\c{c}}ikme{\c{s}}e.
\newblock Real-time quad-rotor path planning using convex optimization and
  compound state-triggered constraints.
\newblock In {\em 2019 IEEE/RSJ International Conference on Intelligent Robots
  and Systems (IROS)}, pages 7666--7673. IEEE, 2019.
\newblock URL: \url{https://doi.org/10.1109/IROS40897.2019.8967706}.

\bibitem{szmuk2020successive}
Michael Szmuk, Taylor~P Reynolds, and Beh{\c{c}}et A{\c{c}}{\i}kme{\c{s}}e.
\newblock Successive convexification for real-time six-degree-of-freedom
  powered descent guidance with state-triggered constraints.
\newblock {\em Journal of Guidance, Control, and Dynamics}, 43(8):1399--1413,
  2020.
\newblock URL: \url{https://doi.org/10.2514/1.G004549}.

\bibitem{reynolds2020dual}
Taylor~P Reynolds, Michael Szmuk, Danylo Malyuta, Mehran Mesbahi, Beh{\c{c}}et
  A{\c{c}}{\i}kme{\c{s}}e, and John~M Carson~III.
\newblock Dual quaternion-based powered descent guidance with state-triggered
  constraints.
\newblock {\em Journal of Guidance, Control, and Dynamics}, 43(9):1584--1599,
  2020.
\newblock URL: \url{https://doi.org/10.2514/1.G004536}.

\bibitem{bock1984multiple}
H~G Bock and K~J Plitt.
\newblock A multiple shooting algorithm for direct solution of optimal control
  problems.
\newblock {\em IFAC proceedings volumes}, 17(2):1603--1608, July 1984.
\newblock URL: \url{https://doi.org/10.1016/s1474-6670(17)61205-9}.

\bibitem{quirynen2015multiple}
Rien Quirynen, Milan Vukov, and Moritz Diehl.
\newblock Multiple shooting in a microsecond.
\newblock In {\em Multiple Shooting and Time Domain Decomposition Methods:
  MuS-TDD, Heidelberg, May 6-8, 2013}, pages 183--201. Springer, 2015.
\newblock URL: \url{https://doi.org/10.1007/978-3-319-23321-5_7}.

\bibitem{nocedal2006numerical}
Jorge Nocedal and Stephen Wright.
\newblock {\em Numerical Optimization}.
\newblock Springer, 2 edition, July 2006.
\newblock URL: \url{https://doi.org/10.1007/978-0-387-40065-5}.

\bibitem{drusvyatskiy2019efficiency}
Dmitriy Drusvyatskiy and Courtney Paquette.
\newblock Efficiency of minimizing compositions of convex functions and smooth
  maps.
\newblock {\em Mathematical Programming}, 178:503--558, 2019.
\newblock URL: \url{https://doi.org/10.1007/s10107-018-1311-3}.

\bibitem{diamond2016cvxpy}
Steven Diamond and Stephen Boyd.
\newblock {CVXPY}: A python-embedded modeling language for convex optimization.
\newblock {\em Journal of Machine Learning Research}, 17(83):1--5, 2016.
\newblock URL: \url{https://www.cvxpy.org/}.

\bibitem{domahidi2013ecos}
Alexander Domahidi, Eric Chu, and Stephen Boyd.
\newblock {ECOS}: An {SOCP} solver for embedded systems.
\newblock In {\em 2013 European control conference {(ECC)}}, pages 3071--3076.
  IEEE, 2013.
\newblock URL: \url{https://doi.org/10.23919/ECC.2013.6669541}.

\bibitem{aps2019mosek}
Mosek ApS.
\newblock Mosek optimization toolbox for matlab.
\newblock {\em User's Guide and Reference Manual, Version}, 4(1), 2019.
\newblock URL: \url{https://docs.mosek.com/10.1/toolbox.pdf}.

\bibitem{yu2022proportional}
Yue Yu, Purnanand Elango, Ufuk Topcu, and Beh{\c{c}}et A{\c{c}}{\i}kme{\c{s}}e.
\newblock Proportional--integral projected gradient method for conic
  optimization.
\newblock {\em Automatica}, 142:110359, 2022.
\newblock URL: \url{https://doi.org/10.1016/j.automatica.2022.110359}.

\bibitem{yu2022extrapolated}
Yue Yu, Purnanand Elango, Beh{\c{c}}et A{\c{c}}{\i}kme{\c{s}}e, and Ufuk Topcu.
\newblock Extrapolated proportional-integral projected gradient method for
  conic optimization.
\newblock {\em IEEE Control Systems Letters}, 7:73--78, 2022.
\newblock URL: \url{https://doi.org/10.1109/LCSYS.2022.3186647}.

\bibitem{cartis2011evaluation}
Coralia Cartis, Nicholas~IM Gould, and Philippe~L Toint.
\newblock On the evaluation complexity of composite function minimization with
  applications to nonconvex nonlinear programming.
\newblock {\em SIAM Journal on Optimization}, 21(4):1721--1739, 2011.
\newblock URL: \url{https://doi.org/10.1137/11082381X}.

\bibitem{chari2024fast}
Govind~M Chari, Abhinav~G Kamath, Purnanand Elango, and Beh{\c{c}}et
  A{\c{c}}{\i}kme{\c{s}}e.
\newblock Fast monte carlo analysis for 6-{DoF} powered-descent guidance via
  {GPU}-accelerated sequential convex programming.
\newblock In {\em AIAA SciTech 2024 Forum}. AIAA, 2024.
\newblock URL: \url{https://doi.org/10.2514/6.2024-1762}.

\bibitem{nmpc2024}
Samet Uzun, Purnanand Elango, Abhinav~G. Kamath, Taewan Kim, and Beh{\c{c}}et
  A{\c{c}}{\i}kme{\c{s}}e.
\newblock Successive convexification for nonlinear model predictive control
  with continuous-time constraint satisfaction.
\newblock In {\em 8th IFAC Conference on Nonlinear Model Predictive Control},
  Kyoto, Japan, August 2024.
\newblock URL: \url{https://doi.org/10.48550/arXiv.2405.00061}.

\bibitem{aircraft2024}
Taewan Kim, Abhinav~G. Kamath, Niyousha Rahimi, Jasper Corleis, Beh{\c{c}}et
  A{\c{c}}{\i}kme{\c{s}}e, and Mehran Mesbahi.
\newblock Approach and landing trajectory optimization for a 6-dof aircraft
  with a runway alignment constraint, 2024.

\bibitem{kamath2023seco}
Abhinav~G. Kamath, Purnanand Elango, Yue Yu, Skye Mceowen, Govind~M. Chari,
  John~M. Carson, III, and Beh{\c{c}}et A{\c{c}}{\i}kme{\c{s}}e.
\newblock Real-time sequential conic optimization for multi-phase rocket
  landing guidance.
\newblock {\em IFAC-PapersOnLine}, 56(2):3118--3125, 2023.
\newblock URL: \url{https://doi.org/10.1016/j.ifacol.2023.10.1444}.

\bibitem{starship2021}
Austin DeSisto.
\newblock Starship and its belly flop maneuver, 2021.
\newblock URL:
  \url{https://everydayastronaut.com/starships-belly-flop-maneuver}.

\end{thebibliography}
}
\section*{Appendix}

\appendix

\section{Construction of the specification in quadrotor flight problem} \label{app:qf}

We first define ${}^f s^1_k$ and ${}^g s^1_k$ functions to check whatever the speed is less than $v_{\mathrm{save}}$ and the quadrotor is in the battery charging station, respectively, from time step $k$ to the time step $k + (k_{\mathrm{w}}-1)$. Note that
\scalebox{0.93}{
$
\begin{aligned}
    f(x) = v_{\mathrm{save}} - \|v\|, \; & \; \varphi_1 := (f(x) \geq 0), \; \vartheta_1 = \bm{G}_{\intv{0}{k_{\mathrm{w}}-1}} \varphi_1\\
    g(x) = d_{\mathrm{w}} - \|r - r_w\|, \; & \; \varphi_2 := (g(x) \geq 0), \; \vartheta_2 = \bm{G}_{\intv{0}{k_{\mathrm{w}}-1}} \varphi_2
\end{aligned}
$}

Then
\begin{align*}
    {}^f s^1_k(x) 
    :=& \Gamma^{ \vartheta_1 }_{\bm{\epsilon}, \bm{p}, \bm{w}} (x, k) \\
    =& {}^{\wedge} h_{p_1, w_1}^{\epsilon_1} ((f(x_k), f(x_{k+1}), \dots, f(x_{k + k_{\mathrm{w}}-1 }))) \\
    {}^g s^1_k(x) 
    :=& \Gamma^{ \vartheta_2 }_{\bm{\epsilon}, \bm{p}, \bm{w}} (x, k) \\
    =& {}^{\wedge} h_{p_2, w_2}^{\epsilon_2} ((g(x_k), g(x_{k+1}), \dots, g(x_{k + k_{\mathrm{w}}-1 })))
\end{align*}
$\forall k \in \intv{1}{K - (k_{\mathrm{w}}-1)}$.

We now define $s_k^2$ function such that the output is positive if and only if speed is less than $v_{\mathrm{save}}$ until the time step $k + (k_{\mathrm{w}}-1)$ as follows:
\begin{align*}
    s_k^2(x) &:= {}^{\wedge} h_{p_3, w_3}^{\epsilon_3} (({}^f s^1_1(x), {}^f s^1_2(x), \dots, {}^f s^1_k(x) ))
\end{align*}
$\forall k \in \intv{1}{K - (k_{\mathrm{w}}-1)}$.

We then define $s^3_k$ function such that the output is positive if and only if speed is less than $v_{\mathrm{save}}$ until the time step $k + (k_{\mathrm{w}}-1)$ and the being in a battery charging station specification is satisfied between the time steps $k$ and $k + (k_{\mathrm{w}}-1)$ as follows:
\begin{align*}
    s_k^3(x) &:= {}^{\wedge} h_{p_4, w_4}^{\epsilon_4} ((s_k^2(x), {}^g s^1_k(x) ))
\end{align*}
$\forall k \in \intv{1}{K - (k_{\mathrm{w}}-1)}$.

The $\vartheta_1 \bm{U}_{[0 : K -( k_{\mathrm{w}}-1) ]} \vartheta_2$ specification is satisfied if and only if $s_k^3(x) > 0$ for any $k \in \intv{1}{K - (k_{\mathrm{w}}-1)}$; therefore,
\begin{align*} \label{eq:app_qf_stl}
    \Gamma^{ \vartheta_1 \bm{U}_{[0 : K -( k_{\mathrm{w}}-1) ]} \vartheta_2 }_{\bm{\epsilon}, \bm{p}, \bm{w}} (x, 1) = {}^{\vee} h_{p_5, w_5}^{\epsilon_5} ((&s^3_1(x), s^3_2(x), \dots, \\ 
    & s^3_{K - (k_{\mathrm{w}}-1)}(x) ))
\end{align*}
\end{document}